%Eric's version of give.
\outer\def\give#1. {\medbreak
             \noindent{\bf#1. }}                     %introduces proofs etc.
\outer\def\section #1\par{\bigbreak\centerline{\S
     {\bf#1}}\nobreak\smallskip\noindent}
\def\({\left(}
\def\){\right)}

\def\sqr#1#2{{\vcenter{\hrule height.#2pt              %qed
     \hbox{\vrule width.#2pt height#1pt\kern#1pt
     \vrule width.#2pt}
     \hrule height.#2pt}}}
\def\square{\mathchoice\sqr{5.5}4\sqr{5.0}4\sqr{4.8}3\sqr{4.8}3}
\def\qed{\hskip4pt plus1fill\ $\square$\par\medbreak}

%section 1

%section 2

%section 4
\def\cA{{\cal A}}
\def\cB{{\cal B}}

\def\cD{{\cal D}}
\def\cE{{\cal E}}

\def\cG{{\cal G}}
\def\cH{{\cal H}}

%general
 
\def\CC{{\rm\kern.24em\vrule width.02em height1.4ex depth-.05ex\kern-.26em C}}
\def\RR{{\,\rm{\vrule width.02em height1.55ex depth-.07ex\kern-.3165em R}}}

     % complex numbers
\def\C{{\bf C}}
     % C to power n
\def\cp1{{{\bf P}^1}}
         % real numbers
\def\R{{\bf R}}
   % R to power n
          % integers
   % ball of dimension n

\def\Z{{\bf Z}}

                                 % closed n-ball
		 % real part
		 % imaginary part
\def\bar{\overline}              % conjugate
     % smooth of order #1
            % tangent space, roman T
              % backslash
                 % absolute value
 
             % disc
             % disk
		 % interior, in math mode
 
                % partial derivative
		 % boundary (same thing)

%\def\nbd{neighborhood\ }         % neighborhood
%\def\vec#1#2{({#1}_1,\dots,{#1}_{#2})}
				 % vector of #1 of length #2
%\def\inner#1#2{\langle{#1},{#2}\rangle}  % inner product

\input epsf.sty
\magnification\magstep1
%\headline{\rm\hfil February 12, 2004}

\centerline{\bf Real Polynomial Diffeomorphisms with Maximal Entropy:}
\centerline{\bf II. Small Jacobian}
\smallskip
\centerline{Eric Bedford and John Smillie}
\bigskip
\section 0.  Introduction

The problem of understanding the dynamical behavior of  diffeomorphisms has played a central
role in the field of dynamical systems. One way of  approaching this question is to ask about
generic behavior in the space of diffeomorphisms.  Another way to approach it is to ask about
behavior in some specific parametrized family.  The family of diffeomorphisms of $\R^2$
introduced by H\'enon has often played the  role of such a test case. This is a two parameter
family given by the formula 
$$f_{a,b}(x,y)=(a-x^2-by,x)$$ for $b\ne0$. There are regions of parameter space which are well 
understood. If we fix $b$ then for $a\ll0$ the nonwandering set of $f_{a,b}$ is empty. For 
$a\gg0$, it is shown in [DN] that the restriction of $f_{a,b}$ to its nonwandering set is
hyperbolic and topologically conjugate to the full two-shift. Such  diffeomorphisms are called
``horseshoes''. How the dynamics changes between these two extremes  has been the subject of
much investigation.  The case $b=0$ is an  interesting special case. In this case the map
$f_{a,b}$ is not a diffeomorphism; in fact the  dynamical behavior is essentially one
dimensional.  The dynamical  complexity of
$f_{a,0}$ increases monotonically with $a$ (see [MT]).  For other values of $b$ no  such
results are known. In fact [KKY] show that in some respects the behavior should  not be
expected to be monotone. One way of measuring the topological complexity is through the
topological entropy,
$h_{top}(f_{a,b})$. This is a continuous real  valued function of the parameters which takes
on values in the interval $[0,\log 2]$.  The case $a\ll0$ corresponds to  $h_{top}=0$. The
case $a\gg0$ corresponds to 
$h_{top}=\log 2$. In this paper we study the set of parameters $(a,b)$ for which 
$h_{top}(f_{a,b})$ takes on its maximal value. We say that $f_{a,b}$ has maximal entropy if 
$h_{top}(f_{a,b})=\log 2$. We analyze the ``maximal entropy locus'' when the Jacobian parameter
$b$ is small. We show:

\proclaim Theorem 1. For each $b$ with $|b|<.08$ there is a unique 
$a=a_b$ so that
$h_{top}(f_{ab})<\log 2$ for $a<a_b$ and $h_{top}(f_{ab})=\log 2$ for 
$a\ge a_b$.
Further, we have:
\item{(1)} If $a>a_b$, $f_{ab}$ is a hyperbolic horseshoe.
\item{(2)} If $a=a_b$, $f_{ab}$ has a
quadratic tangency between stable and unstable manifolds of fixed 
points.  This tangency
is homoclinic when $b>0$ and heteroclinic when $b<0$.

The next result discusses properties of the function $b\mapsto a_b$ defined in Theorem 1.

\proclaim Theorem 2.  The function $b\mapsto a_b$ is continuous on the interval
$(-.08,.08)$.  It is analytic on the subintervals $(-.08,0)$ and $(0,.08)$ but not
differentiable at $b=0$.  Furthermore, there is a generic
unfolding of the homoclinic tangency at the parameter $(a_b,b)$, i.e., at
the point of  tangency, the stable and
unstable manifolds move past one another with positive speed with respect to  $a$.

The terminology ``generic unfolding'' will be explained in greater detail in \S5.

Part of Theorem 1 follows from a more general analysis of polynomial  diffeomorphisms of
maximal entropy in degree $d\ge 2$, which was carried out in  [BS8] and [BS1]. In  particular
we proved in this more general context that a maximal entropy polynomial diffeomorphism  is
either hyperbolic or has a quadratic tangency between stable and unstable manifolds of 
periodic points. The contribution of this paper is to describe the set of parameter values 
corresponding to these two types of behavior.

Though these results are stated for the diffeomorphisms 
$f_{a,b}:\R^2\to\R^2$ our methods give us very complete information about the corresponding
complex extensions
$f_{a,b}:\C^2\to\C^2$ for maximal entropy parameter values. In fact it  is the analysis of
these complex extensions which allows us to obtain information about  the real H\'enon
diffeomorphisms. In particular we take advantage of theory of  intersections of complex
manifolds to analyze the complex extensions of the of the real stable  and unstable manifolds.

In addition to proving Theorems 1 and 2, a goal of this paper is to develop the technique of
crossed mappings as a method of more general applicablilty in the analysis of families of
polynomial diffeomorphisms of ${\bf C}^2$.  These techniques are explored further in [BS3].

We note that this is not the first time that complex methods have  been used to address
similar questions.   J.H. 
Hubbard and R.~Oberste-Vorth [O] used complex methods to improve the result of
Devaney-Nitecki.  And Forn\ae ss and  Gavosto [FG1,2] have
used complex methods to show that there is a generic unfolding of a  complex tangency for
$f_{a,b}$ for certain parameters $(a,b)$.

\section 1. The Quadratic Horseshoe Locus: First Approximation

We consider mappings of the form
$$f_{a,b}(x,y) = (a-x^2-by,x)\eqno(1.1)$$
with $b\ne0$.  Note that $f$ may also be written in the form
$\chi^{-1}\circ f\circ\chi=(y,y^2-a-bx)$, where we set 
$\chi(x,y)=(-y,-x)$.  We
say that
$f_{a,b}$ is a {\it (complex) horseshoe} if
$f_{a,b}$ is hyperbolic on $J=J(f_{a,b})$ and if $f|J$ is 
topologically conjugate
to the full 2-shift.  If, in addition, $J\subset\R^2$, we say that 
$f_{a,b}$ is a
{\it real horseshoe}.  J.H.\ Hubbard and R.
Oberste-Vorth have obtained estimates
on the (complex) horseshoe locus; see [O] and [MNTU, 
Proposition 7.4.6].  These are summarzied in the following:

\proclaim Theorem 1.1.  If $b\ne0$, and if
$|a|>2(1+|b|)^2$, then $f_{a,b}$ is a (complex) horseshoe.  If, in addition,
$b\in\R$, and $a>0$, then $f_{a,b}$ is a real horseshoe.

Since horseshoes have entropy equal to $\log 2$, the following result 
gives a large
region of parameter space where there are no horseshoes.  This is the region to
the left in Figure~1.

\proclaim Theorem 1.2.  Define $\sigma^-(b)=2-{13\over 8} b-{7\over 
64} b^2$ and
$\sigma^+(b)=2+{7\over 4}b + {5\over 16}b^2$.  If $(a,b)\in\R^2$ 
satisfy $b\ne0$,
$|b|\le1$ and
$a\le\max(\sigma^-(b),\sigma^+(b))$, then the entropy of
$f|\R^2$ is less than
$\log2$.

\give Proof.  First we note that a fixed point of $f_{a,b}$ has the form
$(x,y)$, where
$$x=y=-{1\over 2}\left[b+1 \pm\sqrt{(b+1)^2 + 4a}\right]. \eqno(1.2)$$
Now we recall some results from [BS1].  If $f_{a,b} $ is a quadratic
diffeomorphism of $\R^2$, and if $f|\R^2$ has entropy $\log2$, then $f_{a,b}$
has two fixed points, which must be saddles.  Further, one of these points must
be unstably one-sided, and the unstable eigenvalue of $Df$ at this fixed point
is (strictly) greater than 4.  The other fixed point has a negative eigenvalue
in the unstable direction, and this eigenvalue must be less than $-2$.

The differential is given in $(x,y)$-coordinates as
$$Df_{a,b}=\pmatrix{-2x & -b\cr 1& 0\cr}.$$
The product of the eigenvalues is $b$, so we may write them as $\lambda$ and
$b/\lambda$.  Thus the trace of the differential is
$$-2x=\lambda+b/\lambda.$$
If $|b|\le1$, then $\lambda\mapsto\lambda+b/\lambda$ is strictly increasing in
$\lambda$ for $|\lambda|>1$.  Thus the condition that there is an eigenvalue
greater than 4 gives us the inequality
$$-2x>4+{b\over 4},\eqno(1.3)$$
and the condition that there is an eigenvalue less than $-2$ gives the
inequality
$$-2x<-(2+{b\over2}).\eqno(1.4)$$
(Note that the inequalities (1.3) and (1.4) refer to different fixed 
points and thus
involve different values of $x$.)

Now we substitute expression (1.2) into (1.3) and obtain
$$\eqalign{ b+1 \pm \sqrt{(b+1)^2+4a} & > 4+{b\over 4}\cr
   \sqrt{(b+1)^2+4a} & > 3 - {3\over 4}b\cr
b^2 +2b+1+4a > 9 - {9\over2}b + {9\over16} b^2\cr}$$
which gives $a>\sigma^-(b)$.

Similarly, we substitute (1.2) into (1.4) and find
$$\eqalign{ b+1 \pm \sqrt{(b+1)^2+4a} & < -2 - {b\over 2}\cr
\pm \sqrt{(b+1)^2+4a} & < -3 -{3\over 2}b\cr
b^2 +2b + 1 + 4a & > 9 + 9b + {9\over 4}b^2,\cr}$$
where the last inequality is reversed because the quantities being squared are
negative.  Thus $a>\sigma^+(b)$.  The case that one of these inequalities fails
happens exactly when we have $a\le\max(\sigma^+(b),\sigma^-(b))$, and 
in this case
the entropy is not equal to $\log2$.\qed

\epsfxsize4.0in
\centerline{\epsfbox{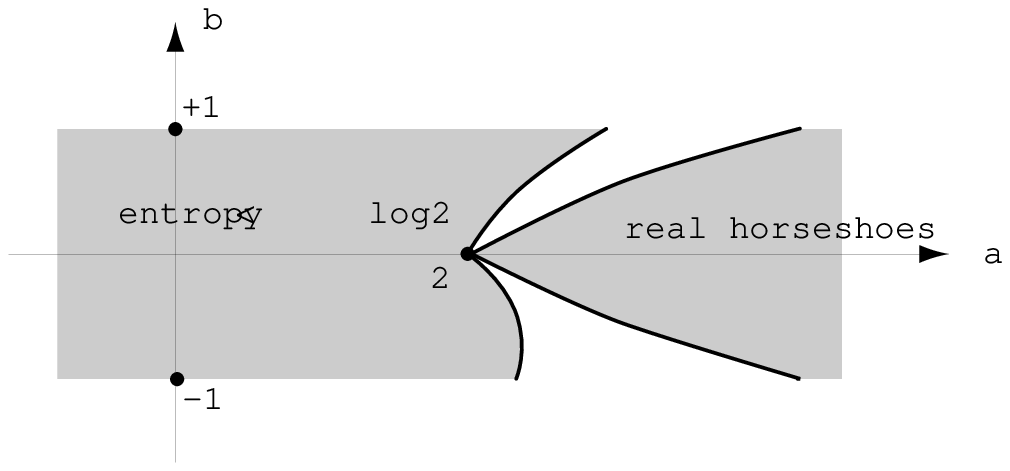}}

\centerline{Figure 1}

Figure 1 shows the information on parameter space that is given by Theorems 1.1 and
1.2.  This figure considers only parameters $|b|\le1$.  In fact,
we restrict our attention without loss of generality to the case
$|b|\le1$ throughout this paper.  Each of the
items discussed in the theorem: maximal entropy, the horseshoe property, and generic
unfolding, will hold for $f$ if and only if it holds for $f^{-1}$.  Thus the fact that
$f_{a,b}^{-1}$ is conjugate to $f_{{a\over b^2}, {1\over  b}}$ means that the regions
that define these dynamical behaviors are invariant under the involution
$(a,b)\mapsto (ab^{-2},b^{-1})$.  In particular, this gives versions of these Theorems
corresponding to the case $|b|>(.08)^{-1}$. 

\section 2.  Complex Boxes and Crossed Mappings

In order to study the system $f|K:K\to K$, we will introduce an open cover by
``boxes'' $B_i$ and study a family of ``crossed mappings'' 
$f_{i,j}:B_i\to B_j$.
We start by working with $p(z)=2-z^2$ and a covering of its Julia set $[-2,2]$.
The Green function for $[-2,2]$ is
$$G(z)=\log\left|{z+\sqrt{z^2-4}\over 2}\right|.$$
For $\lambda>0$, $p$ induces a 2-fold branched covering
$p:\{G<\lambda\}\to\{G<2\lambda\}$. The level sets $\{G=\lambda\}$ are ellipses
with foci at $\pm2$, and the gradient lines (the orthogonal 
trajectories) are given
by the family of hyperbolae with foci at $\pm2$.

Let us fix $c={1\over2}(\sqrt{17}-1)/2\sim 1.5615 $ and
$d={1\over2}(1+\sqrt{7+2\sqrt{17}})\sim2.4523
$.  Let
$E\subset\C$ be the domain bounded by the ellipse with foci at $\pm2$ 
and passing
through $\pm d$.  It follows that
$p(E)$ is the ellipse with foci at $\pm2$ and passing through 
$\pm(d^2-2)$.  We set
$D_0:=\{\zeta\in\C^2:\Re(\zeta)<0\}\cap E$ and
$D_2:=\{\zeta\in\C^2:\Re(\zeta)>0\}\cap E$ as in Figure 2.  It follows that
$p(D_0)=p(D_2)=p(E)-[2,d^2-2)$.  Let
$D_1$ denote the region in $E$ lying between the hyperbolae with foci 
at $\pm2$ and
which pass through $\pm c$ as in Figure 3.  Thus $p(D_1)$ is the region of the
ellipse
$p(E)$ to the right of the hyperbola with foci at $\pm2$ and passing through
$2-c^2$.   We have the following inclusions:
$$D_0\cup D_1\subset p(D_0)=p(D_2), \quad D_2\subset p(D_1).$$
We may also compute certain distances related to these inclusions:
$$\eqalign{ dist(\partial p(D_0),D_0) & = d^2-2-d\cr
dist(\partial p(D_0),D_1) & = 2-c\cr
dist(\partial p(D_1),D_2) &=c^2-2;\cr}\eqno(2.1)$$
While calculating the distances between ellipses can be difficult in general, these
calculations are straightforward because the relevant ellipses are in a confocal
family.  Thus the minimal distances between these ellipses are realized by points on the
real axis.  By the choices of
$c$ and $d$, we have
$$\Delta:=d^2-d-2=2-c=c^2-2\sim.4384.\eqno(2.2)$$

\epsfxsize4.0in
\centerline{\epsfbox{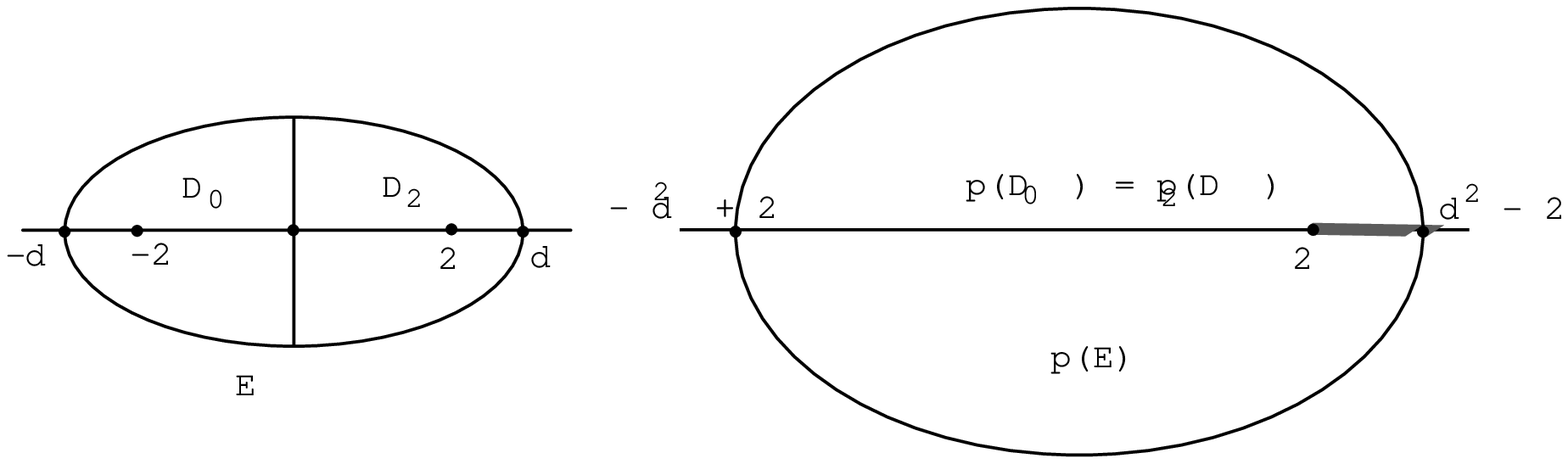}}

\centerline{Figure 2}

\epsfxsize4.0in
\centerline{\epsfbox{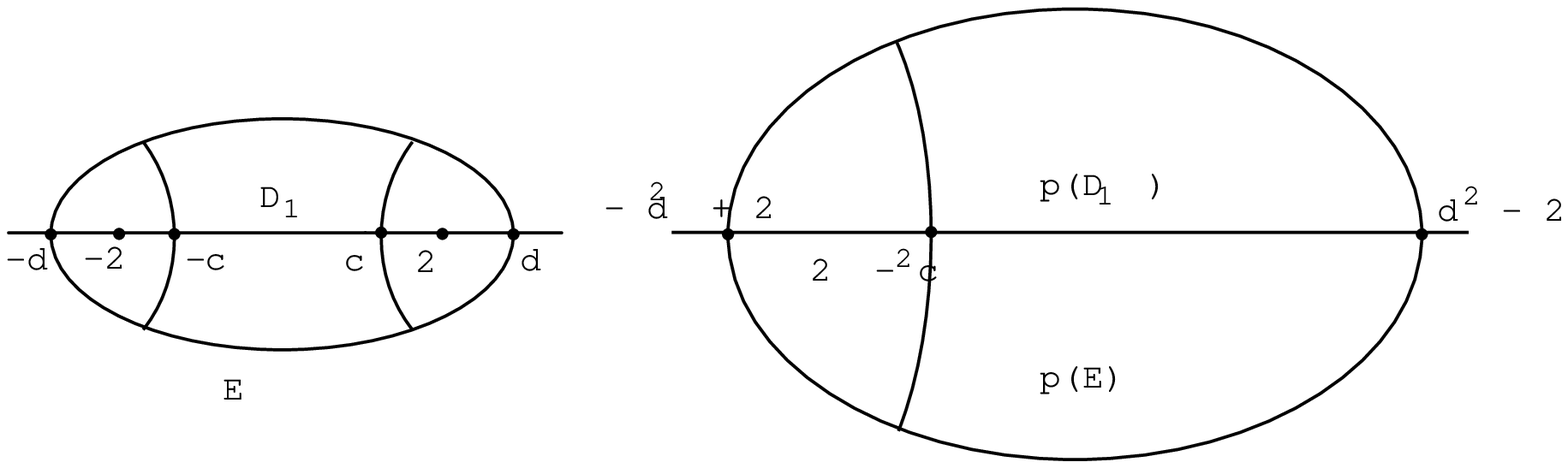}}

\centerline{Figure 3}

Now choose $e>d$ and set $B_j=\{(x,y)\in\C^2:x\in D_j, 
|y|<e\}=D_j\times\{|y|<e\}$
for $j=0,1,2$.  Thus $ B_0\cup B_1\cup B_2 = E\times\{|y|<e\}$.  We 
introduce the set
$$\cA:=\{a,b\in\C,\ \ b\ne0, \ \ |a-2|+e|b|<\Delta\}\approx \{|a-2|+2.4|b|<.43\}. 
\eqno(\ddag)$$

\proclaim Proposition 2.1.   If $(a,b)\in\cA$, then $K\subset 
B_0\cup B_1\cup B_2$.

\give Proof.  By [MNTU, p.\ 238], we know that $K$ is contained in the bidisk
$\{|x|,|y|<R\}$, where $R$ is the larger root of the equation 
$t^2-(1+|b|)t-|a|=0$.
By the condition $|a-2|+e|b|<\Delta$, we conclude that we may take 
$e$ sufficiently
close to $d$, then we have
$$R\le {1+\Delta/e +\sqrt{(1+\Delta/e)^2+4(2+\Delta)}\over2}\sim 2.25845$$
Recall that $pE$ is an ellipse with foci at $\pm2$ and major axis of
length $d^2-2\sim 4.01378$.  We then compute that its minor axis has length
$\sqrt{(d^2-2)^2-4}\sim3.48$.

To prove the Proposition, we need to show that if $(x,y)\in\{|x|,|y|<R\}$, and
if $x\notin E$, then $(x,y)\notin K$.   For such $(x,y)$, the $x$-coordinate of
$f(x,y)$ satisfies:
$$|\pi_v(f(x,y))-p(x)|< |a-2| +|by|<|a-2| +R|b|<\Delta 
+(R-e)|b|<1.03465\Delta$$
since $|b|<\Delta/e$.  Now let $D:=\{\zeta\in pE:dist(\zeta,\partial
(pE))>1.03465\Delta\}$.  Since $x\notin E$, it follows that $px\notin 
pE$, and so
the
$x$-coordinate of $f(x,y)$ does not belong to $D$.  On the other 
hand,  the minor
axis of
$pE$ is 3.48, so that $D$ contains the disk of radius 
$3.48-1.03465\Delta\sim3.0264
>R$.  Thus $f(x,y)\notin K$.
\qed

The vertical and horizontal
components of the boundaries are defined to be
$$\partial_vB_j=(\partial
D_j)\times\{|y|\le e\},\quad\partial_hB_j=\bar D_j\times\{|y|=e\}.$$
We set $\cG=\{(0,0),(0,1), (1,2), (2,1), (2,0)\}$, and we interpret $\cG$ as a
graph on the vertices $\{B_0,B_1,B_2\}$ as in Figure 4.

\epsfxsize4.0in
\centerline{\epsfbox{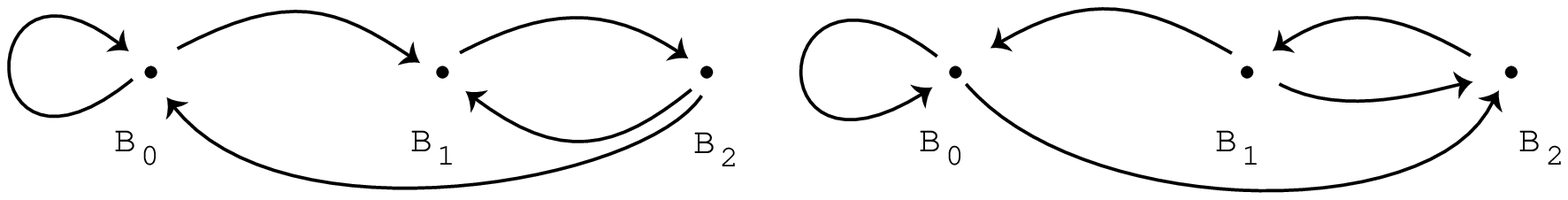}}

\centerline{Figure 4: Graphs $\cG$ for $f$ (on left) and $\cG^{-1}$ 
for $f^{-1}$
(on right)}
\proclaim Proposition 2.2.   If $(\ddag)$ holds, then 
$f_{a,b}(\partial_vB_i)\cap
\bar B_j=\emptyset$ and $f_{a,b}(\bar B_i)\cap \partial_hB_j=\emptyset$ for all
$(i,j)\in\cG$. 

\give Proof.  By estimates (2.1) and (2.2) and the fact that $p(\partial
D_0)=p(\partial D_2)$, we have
$$dist(p(\zeta),\partial D_j)\ge\Delta$$
for $\zeta\in\partial D_i$ and $(i,j)\in\Gamma$.  Thus if 
$x\in\partial D_i$ and
$|y|<e$, the first coordinate of
$f_{a,b}(x,y)$ satisfies
$$|a-x^2-by-p(x)|\le|a-2|+|by|<\Delta.$$
This gives $a-x^2-by\notin\bar D_j$, so
$f(\partial_vB_i)\cap\bar B_j=\emptyset$.

Note that $\partial_hB_j\subset\{|y|=e\}$.  The second coordinate
of $f$ is $x$, so the second condition follows from the fact that $\bar
D_j\cap\{|y|=e\}=\emptyset$, independently of $a$ and $b$.
\qed

Let $\pi_v(x,y)=x$ and $\pi_h(x,y)=y$ denote the projections in the 
vertical and
horizontal directions.  We let $f_{i,j}$ denote the mapping
$f:B_i\cap f^{-1}B_j\to B_j$.  Following [HO], we say that $f_{i,j}$ is a
{\it crossed mapping} if for each $y\in\{|y|<e\}$,
$$\pi_v\circ f:(D_i\times\{y\})\cap f^{-1}B_j\to D_j\eqno(2.4)$$
is proper.  Given $(i,j)\in\cG$, then it follows from Proposition 2.2 
that $f_{ij}$
is a crossed mapping.  We say that the degree of
$f_{ij}$ as a crossed mapping is the mapping degree of the map in 
(2.4) (which is
independent of
$y$).  Similarly, we say that
$f^{-1}:B_j\cap f B_i\to B_i$, which we denote by $f^{-1}_{ji}$, is a crossed
mapping if for each
$x\in D_j$,
$$\pi_h\circ f^{-1}:(\{x\}\times\{|y|<e\})\cap fB_i\to \{|y|<e\}\eqno(2.5)$$
is proper.  As was observed in [HO], $f_{ij}$ is a crossed mapping 
if and only if
$f_{ji}^{-1}$ is.  And
the degree of
$f^{-1}_{ji}$ as a crossed mapping is defined as the mapping degree 
of the map in
(2.5) (which is independent of $x$).  This, in turn, is the same as 
the degree of
$f_{ij}$.  We will say that $(\cB,\cG)$ is a system of crossed mappings, if $f_{i,j}$
induces a crossed mapping from $B_i$ to  $B_j$ for each $(i,j)\in\cG$.   The
following Corollary is a consequence of Proposition 2.2.

\proclaim Corollary 2.3.   If $(\ddag)$ holds, then $(\cB,\cG)$ is a system of crossed mappings.

We define an {\it orbit} in a system of crossed mappings as a bi-infinite sequence
$(p_j,i_j)_{j\in\Z}$ such that for all $j\in\Z$, $p_j\in B_{i_j}$, $(i_j,i_{j+1})\in\cG$, and
$f(p_j)=p_{j+1}$.  Next we give conditions that guarantee that every $f$-orbit
$(p_j)_{j\in\Z}$ in $K$ can be lifted to an orbit of the system of crossed mappings.

\proclaim Proposition 2.4.  Suppose that $K\cap (B_0\cup B_1)\subset f(B_0\cup
B_2)$ and $K\cap B_2\subset f(B_1)$.  Then for $q\in K$ there is an admissible
sequence $I=(i_n)_{n\in\Z}$ such that $f^nq\in B_{i_n}$ for all $n\in\Z$.

\give Proof.  Let us start by making a sequence $J_M=\{j_n: -M\le n\le M\}$ of
finite length.  If we have determined $j_n$ already, then $f^n(q)\in 
B_{j_n}\cap
K$.  If
$j_n=0$ or 1, then by hypothesis $f^{n-1}q\in (B_0\cup B_2)\cap K$. 
Thus we may
choose
$j_{n-1}\in\{0,2\}$ such that $f^{n-1}q\in B_{j_{n-1}}$, and in 
either case we have
$(j_{n-1},j_n)\in\cG$.  Similarly, if $j_n=2$, then $f^nq\in B_2\cap K$, and by
hypothesis we have $f^{n-1}q\in B_1$.  Thus we set $j_{n-1}=1$ and
$(j_{n-1},j_n)=(1,2)\in\cG$.  Starting at $n=M$, we continue 
backwards and generate
an admissible sequence $J_M$.

Now we have admissible sequences $J_1,J_2,\dots$ of increasing 
length.  For each
$M$, there is a sequence $I_M$ that is a subsequence of infinitely 
many sequences
$J_{k_m}$.  We may make $M$ increasingly large and thus obtain an 
infinite sequence
$I$.\qed

\proclaim Proposition 2.5.  If $a,b\in\R$ and if (\ddag) holds, then
$$(\bar B_{0,r}\cup\bar B_{1,r})\cap f(B_0\cup B_1\cup B_2)\subset 
f(B_0\cup B_2)$$
$$\bar B_{2,r}\cap  f(B_0\cup B_1\cup B_2)\subset f(B_1).$$

\give Proof.  We note that $(B_0\cup B_1\cup B_2 -B_0\cup B_2)\cap\R^2 =
\{0\}\times (-e,e)$.  Thus to prove the first inclusion, it suffices 
to show that
$(\bar B_{0,r}\cup\bar B_{1,r})\cap f(\{0\}\times(-e,e))=\emptyset$.  But $\bar
B_{0,r}\cup\bar B_{1,r}=[-d,c]\times[-e,e]$, and the $x$-projection of the
$f$-image of this set is
$$\pi_v\circ f(\{0\}\times(-e,e))=\{a-x^2-by:x=0,|y|<e\}\subset
(2-\Delta,2+\Delta).$$
On the other hand, $\bar B_{0,r}\cup\bar B_{1,r}=[-d,c]\times[-e,e]$. 
Thus $(\bar
B_{0,r}\cup\bar B_{1,r})\cap f(\{0\}\times(-e,e))=\emptyset$ since 
$c+\Delta=2$,
which proves the first inclusion.

For the second inclusion, we note that
$$(B_0\cup B_1\cup B_2-B_1)\cap\R^2=((-d,-c)\cup(c,d))\times(-e,e).$$
The $x$-projection of the $f$-image of this set is
$$\{a-x^2-by:c<|x|<d,|y|<e\}\subset(2-d^2-\Delta,2-c^2+\Delta).$$
On the other hand, the $x$-projection of $\bar B_{2,r}$ is $[c,d]$, which is
disjoint from $(-\infty,2-c^2+\Delta)$ since $2-c^2+\Delta=c$.\qed

Let $V\subset B_i$ be a complex subvariety.  We say that $V$ is a 
{\it horizontal
multi-disk} (resp.\ {\it vertical multi-disk}) if each component of $V$ is
conformally equivalent to a complex disk, and if 
$\partial_hB_i\cap\bar V=\emptyset$
(resp.\ $\partial_vB_i\cap\bar V=\emptyset$).
{With this terminology the union of horizontal (resp.\ vertical) multi-disks is
again a horizontal (resp.\ vertical) multi-disk.}   We denote the set 
of horizontal
(resp.\ vertical) multi-disks by
$\cD_h(B_i)$ (resp.\
$\cD_v(B_i)$).  If $V\in\cD_h(B_0)$ (resp.\ $V\in\cD_v(B_i)$), then 
$\pi_v:V\to D_i$
(resp.\
$\pi_h:V\to\{|y|<e\}$) is proper, and we let
$\delta(V)$ denote the degree of the corresponding projection.  By 
$\cD_h^{m}(B_i)$
(resp.\
$\cD_v^{m}(B_i)$) we denote the set of horizontal (resp.\ vertical) 
multi-disks $V$
such that for each component $W$ of $V$, the degree $\delta(W)$ is no
greater than $m$.  We note the following:
$${\rm If \ }V'\in\cD_{v}(B_i){\rm\  and\ } V''\in\cD_h(B_i), {\rm\ then\ }
\#(V'\cap V'')=\delta(V')\delta(V'').\eqno(2.6)$$
If $f_{ij}$ is a crossed map, and if $V\subset B_i$ is a subvariety for which
$\partial_hB_i\cap\bar V=\emptyset$, then $\tilde f_{ij}(V):=f(V)\cap 
B_j$ is closed
in
$B_j$ and satisfies $\partial_hB_j\cap\bar{\tilde V}=\emptyset$, and is thus a
horizontal subvariety.  If $\deg (f_{ij})$ denotes the degree as a 
crossed map, then
we have
$$\deg(f_{ij})\delta(V)=\delta(\tilde f_{ij}(V)).\eqno(2.7)$$

\proclaim Proposition 2.6.  If $(\ddag)$ holds, it follows that
$$\tilde
f_{12}:\cD_h^m(B_i)\to\cD_h^{2m}(B_j){\rm\ and\ }\tilde
f^{-1}_{21}:\cD_v^m(B_j)\to\cD_v^{2m}(B_i),$$
and if $(i,j)\in\cG$, $(i,j)\ne(1,2)$, then
$$\tilde
f_{ij}:\cD_h^m(B_i)\to\cD_h^m(B_j){\rm\ and\ }\tilde
f^{-1}_{ji}:\cD_v^m(B_j)\to\cD_v^m(B_i).$$

\give Proof.  We will show that $\tilde f_{ij}(V)$ is conformally 
equivalent to a
disk; the degree is given by (2.7).  Suppose
$V$ is a horizontal disk in
$B_i$.  Then, taking boundary inside $\C^2$, we have $\partial
V\subset\partial_v(B_i)$.  By Proposition 2.2,
$f(\partial V)\cap \partial B_j=\emptyset$.  Thus each component of 
$f(V)\cap B_j$
is closed in $B_j$.  The second part of Proposition 2.2 implies that 
$\pi_v:f(V)\cap
B_j\to D_j$ is proper.  Finally, we need to show that each component $W$ of
$f(V)\cap B_j$ is conformally equivalent to the disk.  Since $V$ is a 
disk, there
is a conformal equivalence $\varphi:\Delta\to V\subset\C^2$.  Now the 
components of
$fV\cap B_j$ correspond to the components of
$\{\zeta\in\Delta:f\circ\varphi(\zeta)\in B_j\}=\{\zeta\in\Delta:\pi_v\circ
f\circ\varphi(\zeta)\in D_j\}$.  Since $\bar D_j$ is simply 
connected, there is a
subharmonic function $s$ on $\C$ such that $\bar D_j=\{s\le0\}$.  It 
follows from
the maximum principle that each component of $(\pi_v\circ 
f\circ\varphi)^{-1}(\bar
D_j)=\{s\circ f\circ\varphi\le0\}$ is simply connected.  Finally, 
since $\pi_v\circ
f\circ\varphi$ is an open mapping, each component of $(\pi_v\circ
f\circ\varphi)^{-1}D_j$ is simply connected.
\qed

Since $f_{00}$ is a crossed mapping from $B_0$ to itself, there is a 
saddle fixed
point $p_0\in B_0$.  Let us define $W^{s/u}_0$ to be the connected component of
$W^{s/u}(p_0)\cap B_0$ which contains $p_0$.  It follows that
$W^{s/u}_0\in\cD^1_{v/h}(B_0)$.  Note that
$$W^u_0=\bigcap_{n\ge0}f^nB_0,\quad W^s_0=\bigcap_{n\ge0}f^{-n}B_0.\eqno(2.8)$$
There is also a saddle point $p_1\in B_1\cap B_2$.  We let $W^u_1$ denote the
component of of $W^u(p_1)\cap B_1$ that contains $p_1$.  We show in Proposition 4.3 that if
(\ddag) holds, then it is a horizontal disk of degree 1.

Let us say that a sequence $I=i_0i_1\cdots i_n$ is {\it admissible} if
$(i_k,i_{k+1})\in\cG$ for all $k$.  We will sometimes also say that a sequence
$J=j_0j_1\cdots j_m$ is admissible if $(j_k,j_{k+1})\in\cG^{-1}$ for 
all $k$.  It
will be clear from context whether we mean $\cG$ or $\cG^{-1}$.  For admissible
sequences $I$ (for $\cG$) and $J$ (for $\cG^{-1}$), we use the notation
$$W^u_I=W^u_{i_0i_1\cdots i_n}=\tilde f_{i_{n-1}i_n}\tilde 
f_{i_{n-2}i_{n-1}}\cdots
\tilde f_{i_0i_1}(W^u_{i_0})\eqno(2.9)$$
$$W^s_J=W^s_{j_0j_1\cdots j_n}=\tilde f^{-1}_{j_{n-1}j_n}\tilde
f^{-1}_{j_{n-2}j_{n-1}}\cdots
\tilde f^{-1}_{j_0j_1}(W^s_{j_0}).\eqno(2.10)$$
It follows that $W^{s/u}_0$ are vertical/horizontal disks of degree 1 in $B_0$, and $W^s_{02}$
is a vertical disk of degree 1 in $B_2$.  By Proposition 2.6,
$W^u_{01}$ are vertical/horizontal disks of degree 1; and $W^u_{012}$ is a
horizontal disk of degree 2.  This last statement includes two possibilities:
$W^u_{012}$ might consist of two disjoint disks of degree 1 or one disk on which
$\pi_v$ has degree 2.  In either case, $W^u_{012}$ intersects $W^s_{02}$ in
$B_2$ with multiplicity two, which means that either $W^u_{012}\cap W^s_{02}$
consists of two distinct points, or the intersection is tangential.

\section 3.  Mappings of Real Boxes

Here we work under the additional condition that $f_{a,b}$ is a real 
mapping.  In this section, we will restrict our attention to the real parameter region
$$\cA_r:=\cA\cap\R^2.$$
Let
$\tau$ be the involution of $\C^2$ defined by
$\tau(x,y)=(\bar x,\bar y)$.  The fixed point set of $\tau$ is $\R^2$.  The
condition that $a,b\in\R$ is equivalent to the condition that 
$f_{a,b}$ commutes
with $\tau$.  We say that a set $S\subset\C^2$ is {\it real} if $\tau 
S=S$.  For instance $\tau
B_i=B_i$, so in this terminology $B_i$ is real.  Let
$\cD_{h/v,r}(B_i)$ denote the set of horizontal/vertical disks in 
$\cD_{h/v}(B_i)$
which are real.  If $(a,b)\in\cA_r$, then Proposition 2.6 applies 
to real disks
to yield
$$\tilde
f_{12}:\cD_{h,r}^m(B_i)\to\cD_{h,r}^{2m}(B_j){\rm\ and\ }\tilde
f^{-1}_{21}:\cD_{v,r}^m(B_j)\to\cD_{v,r}^{2m}(B_i),$$
and
$$\tilde
f_{ij}:\cD_{h,r}^m(B_i)\to\cD_{h,r}^m(B_j){\rm\ and\ }\tilde
f^{-1}_{ji}:\cD_{v,r}^m(B_j)\to\cD_{v,r}^m(B_i)$$
for $(i,j)\in\cG$, $(i,j)\ne(1,2)$.

We set $B_i^r:=B_i\cap\R^2$, which
is a rectangle in $\R^2$ with sides parallel to the axes.
\proclaim Proposition 3.1.  If $V\in\cD_{h,r}(B_i)$, then
$V\cap B_{i,r}$ consists of a nonempty, connected, one-dimensional 
curve.  In fact,
there is a conformal uniformization $h:\Delta\to V$ such that
$h(\bar\zeta)=\tau\circ h(\zeta)$.

\give Proof.  Let $\varphi:\Delta\to V$ be a conformal uniformization 
of $V$.  It
follows that
$\kappa:\Delta\ni\zeta\mapsto\varphi^{-1}\circ\tau\circ\varphi(\zeta)\in\Delta$ 
is
an anti-conformal involution of $\Delta$. It follows that $\kappa$ is an
orientation-reversing isometry for the Poincar\'e metric, so the 
fixed point set
$\gamma:=\{\zeta\in\Delta:\kappa(\zeta)=\zeta\}$ is a Poincar\'e geodesic.  Let
$\psi$ be a conformal automorphism of $\Delta$ which maps the real axis
$(-1,1)\subset\Delta$ to $\gamma$.  It follows that 
$\psi^{-1}\circ\kappa\circ\psi$
is an isometric involution of $\Delta$ which fixes $(-1,1)$, so it is 
simply the
map $\zeta\mapsto\bar\zeta$.  Thus
$h=\varphi\circ\psi$ is the desired uniformization.
\qed
If $f$ is a real map,
then for $(i,j)\in\cG$, $f_{ij}$ is a crossed mapping of the pair
$(B_i^r,B_j^r)$.

\proclaim Proposition 3.2.  If $a,b\in\cA_r$, then $B_{0,r}\cap
fB_0$ lies below $B_{0,r}\cap fB_2$ inside $B_{0,r}$, and 
$B_{1,r}\cap fB_0$ lies
below
$B_{1,r}\cap fB_2$ inside $B_{1,r}$.  In particular, let 
$I=0i_1\cdots i_n00$ and
$J=0j_1\cdots j_m20$ be admissible sequences.  Then $W^u_I$ lies below $W^u_J$
inside $B_{0,r}$.  Similarly, if $K=0k_1\cdots k_n01$ and 
$L=0l_1\cdots l_m21$ are
admissible sequences, then $W^u_K$ lies below $W^u_L$ inside $B_{1,r}$.

\give Proof.  The $y$-coordinate of $f$ is $\pi_h\circ f=x$.  Since 
$B_0$ lies to
the left of $B_2$, it follows that the $y$-coordinate of $fB_{0,r}$ 
is less than
that of $fB_{2,r}$.  Thus it lies below.

For the assertions about the pieces of unstable manifolds, we note that if
$I$ is a sequence that ends in $ij$, then $W^u_I\subset fB_i\cap 
B_j$.  Thus for a
sequence $I$ which ends in $00$ and a sequence $J$ which ends in 
$20$, we will have
$W^u_I\subset B_{0,r}\cap fB_0$ which lies below $W^u_J\subset 
B_{0,r}\cap fB_2$.
\qed

If $(i,j)\in\cG$, $(i,j)\ne(1,2)$, then the crossed mapping $f_{ij}$ 
has degree 1.
This means that real, horizontal curves in $B_{i,r}$ which run from 
left to right
are taken to real, horizontal curves in $B_{j,r}$ which run either from left to
right or from right to left.  If the left-to-right direction is preserved, we
assign the symbol $\epsilon_u=+$ to $f$.  Otherwise, we set $\epsilon_u=-$.
Similarly, real, vertical curves in $B_{j,r}$ which run from bottom to top are
mapped under
$f^{-1}$ to real, vertical curves which either run from bottom to top 
or from top to
bottom.  If the run from bottom to top, then we assign the symbol 
$\epsilon_s=+$ to
$f$.  Otherwise, $\epsilon_s=-$.

\epsfxsize4.5in
\centerline{\epsfbox{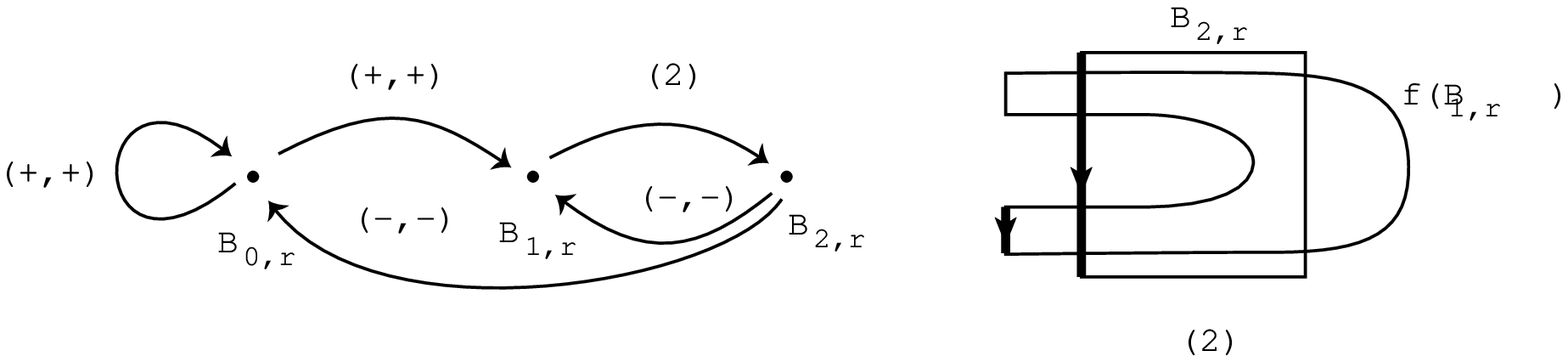}}
\centerline{Figure 5:  Graph induced by $f$ (orientation-preserving)}

\epsfxsize4.5in
\centerline{\epsfbox{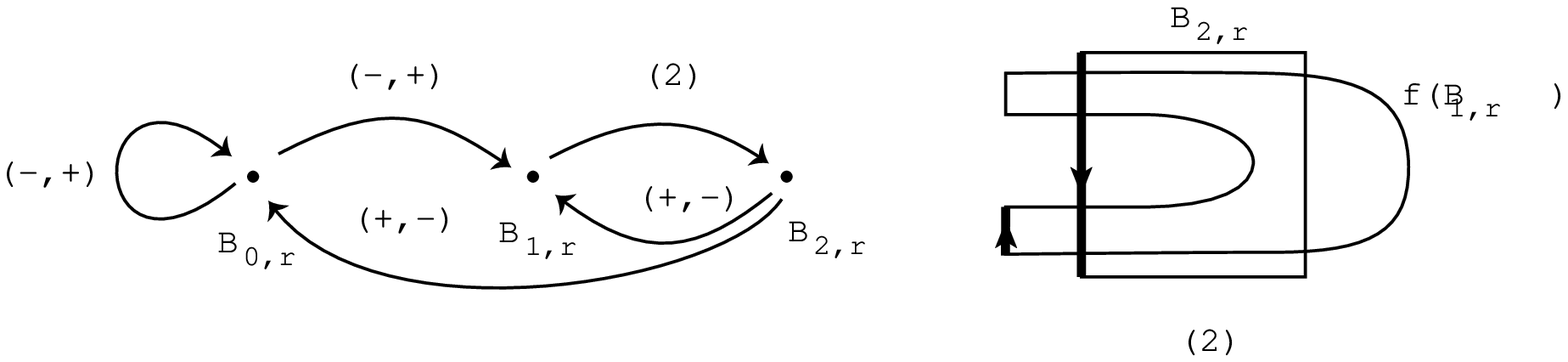}}
\centerline{Figure 6:  Graph induced by $f$ (orientation-reversing)}

\proclaim Proposition 3.3.  If $a,b\in\cA_r$, then the signs
$(\epsilon_s,\epsilon_u)$ are given as in Figures 5 and 6.

\give Proof.  First we consider the degenerate case $b=0$.  The map
$a-x^2=\pi_v\circ f_{a,0}(x,y)$ is increasing on $D_0\cap\R=(-d,0)$ 
and decreasing
on
$D_2\cap\R=(0,d)$.  Thus we have $\epsilon_u=+$ on $D_0\cap\R$ and 
$\epsilon_u=-$
on $D_0\cap\R$.  This condition continues to hold for $b\ne0$.  Thus we have
$\epsilon_u=+$ on $D_0$ and
$\epsilon_u=-$ on $D_2$.  This continues to hold for $b\ne0$, so the arrows of
$\cG$ emanating from $B_0$ should be labeled
$(\cdot,+)$, and the arrows emanating from $B_2$ should be labeled $(\cdot,-)$.
In the orientation-preserving case, the only possible labels are $(+,+)$ and
$(-,-)$.  In the orientation-reversing case, the only possible labels 
are $(+,-)$
and $(-,+)$.  Thus we have the labeling shown in the graphs in Figures 5 and 6. \qed

The crossed map $f_{12}$ has degree 2 and is less easy to work with.  The
illustrations on the right hand sides of Figures 5 and 6 indicate its
combinatorial behavior in the following sense.  The left side of the
vertical boundary of $B_{1,r}$ is
$\{-c\}\times[-e,e]$, and the right side is $\{c\}\times[-e,e]$.  In 
the degenerate
case $b=0$, $f_{a,0}$ maps the left boundary to the point 
$(a-c^2,-c)$ which is to
the left of $B_{2,r}$; and the right boundary goes to $(a-c^2,c)$ 
which is directly
above
$(a-c^2,-c)$.  If $b\ne0$, then the image of the left boundary will 
continue to be
to the left of $B_{2,r}$ and below the image of the right boundary.  The use we
make of this combinatorial/topological information is given in Proposition 3.4,
whose proof is a straightforward consequence of the preceding discussion.

\proclaim Proposition 3.4.   Suppose $a,b\in\cA_r$.  Suppose,
too, that $A_1$ and $A_2$ are curves that cross $B_{1,r}$ from left 
to right and
that $A_1$ lies below $A_2$ inside $B_{1,r}$.  If $f$ preserves 
orientation, then
the curves $C_1=\tilde
f_{12}A_1$ and $C_2=\tilde f_{12}A_2$ open to the left, and $C_2$ 
lies inside $C_1$
as illustrated in Figure 7.  If $f$ reverses orientation, then the relative
positions of $C_1$ and
$C_2$ are exchanged.

\bigskip
\epsfxsize4.5in
\centerline{\epsfbox{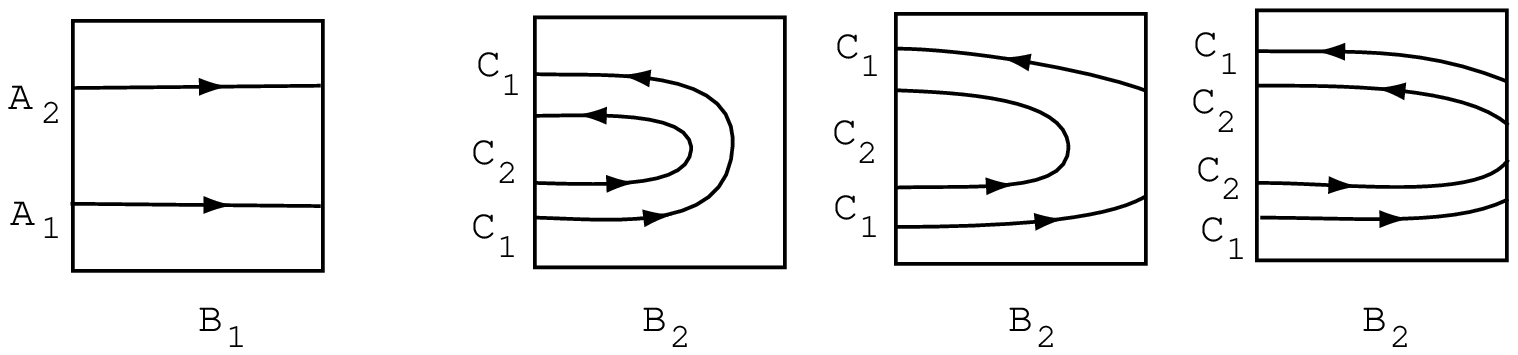}}
\centerline{Figure 7:  Curve $C_2$ lies inside $C_1$: three possibilities.}
\bigskip
In the sequel, we will work with parameter values in $\cA_r$.  However, for many of the
arguments  the essential point is that $(\cB,\cG)$ is a system of crossed mappings with
the ``combinatorial'' behavior given in Figures 5 and 6.  Thus we are led to the
following condition:
$$\eqalign{&(\cB,\cG)\ \hbox{\rm is a family of real, crossed mappings, with the} \cr
&{\rm\
topological\ configurations\ shown\ in\ Figures\ 5,\ 6,\ and\ 7.}}\eqno(\dag)$$
We may summarize the discussion above by the statement: {\sl If $a,b\in\R$ and (\ddag) holds,
then (\dag) holds.}   We also introduce the two conditions
$$\#(W^s_{02}\cap W^u_{012}\cap B_{2,r})=2{\rm \ if\ }b>0, {\rm\ and\ }
\#(W^s_{02}\cap W^u_{12}\cap B_{2,r})=2{\rm\ if\ }b<0.
\eqno(*)$$
$$\#(W^s_{02}\cap W^u_{01212}\cap B_{2,r})=4{\rm \ if\ }b>0, {\rm\ and\ }
\#(W^s_{02}\cap W^u_{12012}\cap B_{2,r})=4 {\rm \ if\ }b<0.\eqno(**)$$
\noindent{\bf Remark on notation.}  We have now defined a parameter domain $\cA_r$ as
well as three conditions that may or may hold for a given parameter value
$(a,b)$.  The condition (\dag) requires
the boxes
$\cB$ to have specified behavior under $f$ and $f^{-1}$.  The conditions $(*)$ and
$(**)$ define dynamical characteristics of $f_{a,b}$.  It will be shown below that
(\dag) holds for all parameters in $\cA_r$ and that $(**)$ implies $(*)$.
\proclaim Proposition 3.5.  If (\dag) holds, then $(**)\Rightarrow(*)$.

\give Proof.  We will treat only the case $b<0$ since the case $b>0$ 
is similar.
Let us suppose that $(*)$ fails.  We map $W^u_1$ forward under
$\tilde f_{12}$ to $W^u_{12}$.  By  Proposition 3.1, $W^u_{12}\cap 
B_{2,r}$ is a
nonempty, connected curve, and by Proposition 3.4  it forms a curve 
which opens to
the left, which by hypothesis does not intersect
$W^s_{02}$.  This is pictured in the pair of boxes on the left hand 
side of Figure 8.  Next we map
$W^u_{12}$ forward under
$\tilde f_{20}$. Again by Proposition 3.1, $W^u_{120}\cap B_{0,r}$ is 
a nonempty,
connected curve.  Since the sign of
$f_{20}$ is
$(\cdot,-)$, the
$x$-direction of the curve is reversed, so $W^u_{120}\cap B_{0,r}$ opens to the
right.  By Proposition 3.2,
$W^u_{120}$ lies above $W^u_0=W^u_{00}$, which is drawn in gray as a 
visual aid to
the reader, although it is not necessary for the proof.  (The gray dot is
$p_0$, and $W^u_1$ is above $W^u_{01}$ in $B_{0,r}$ by Proposition 3.2.)  Since the
sign of $f_{21}$ is $(+,\cdot)$, the vertical orientation is preserved, so
$W^u_{121}$ contains $W^u_1$ as well as a curve below it.  Since the sign of
$f_{20}$ is $(+,\cdot)$, the vertical orientation is preserved, so 
the upper part
of $W^u_{120}$ with a single hash mark is identified with $W^u_1$ on the set
$B_{0,r}\cap B_{1,r}$.  Since
$W^s_{02}\cap W^u_{12}=\emptyset$, $W^u_{120}$ is disjoint from $W^s_0$, so we
obtain the picture as in the right hand pair of boxes in Figure 8.

\bigskip
\epsfxsize4.5in
\centerline{\epsfbox{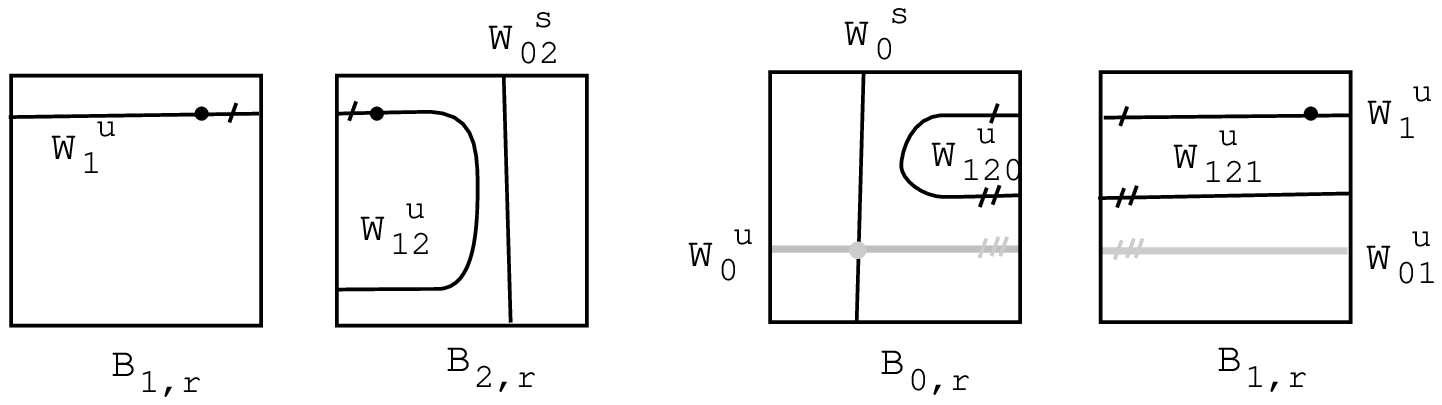}}
\centerline{Figure 8}
\medskip

Next we map $W^u_{120}$ forward under $\tilde f_{01}$.  This is shown 
in the left
hand picture of Figure 9.  Since $f_{01}$ has signature $(-,\cdot)$, 
the vertical
orientation is reversed, so $W^u_{1201}$ lies below $W^u_{121}$ and 
$W^u_{01}$ in
$B_{1,r}$.  Finally, we map forward under
$\tilde f_{12}$ and obtain the picture in the right hand box of 
Figure 9.  The two
arches of $W^u_{12012}$ lie inside $W^u_{12}$ by Proposition 3.4.  Thus
$W^u_{12012}$ cannot intersect $W^s_{02}$, so condition $(**)$ does not hold.
\qed

\bigskip
\epsfxsize2.9in
\centerline{\epsfbox{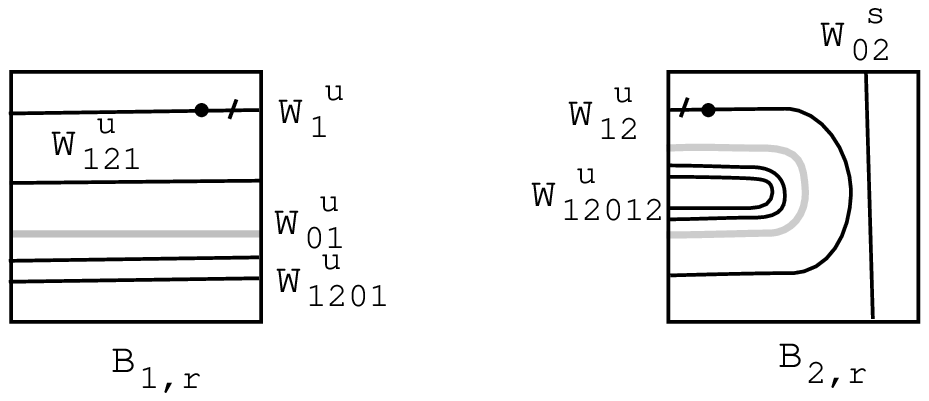}}
\centerline{Figure 9}
\medskip
Figure 10 illustrates conditions $(*)$ and $(**)$ in the case $b>0$.  To understand
Figure 10, start in the left hand box with $W^u_0$ and $W^s_0$ 
passing through the
saddle point $p_0$.  We move $W^u_0$ to box $B_{1,r}$ via the map 
$\tilde f_{01}$,
and to box $B_{2,r}$ via $\tilde f_{02}$. The map $f_{02}$ has degree 2, and
$\tilde f_{02}W^u_{01}=W^u_{012}$ is a curve of degree 2 which opens 
to the left by
Proposition 3.4.  By condition $(*)$, $W^u_{012}$ crosses $W^s_{02}$. 
The crossed
map $f_{20}$ has degree 1 and sign $(\cdot,-)$, so the left-opening, degree two
curve
$W^u_{012}$ produces a degree two curve $W^u_{0120}=\tilde f_{20}W^u_{012}$ in
$B_{0,r}$ which opens to the right.  Condition $(*)$ maps forward under 
$f_{20}$, so
$W^u_{0120}$ intersects $W^u_0$.

The crossed map $f_{21}$ has degree 1, so $\tilde
f_{21}W^u_{012}=W^u_{0121}$ has degree 2 and by Proposition 3.2, it  lies above
$W^u_{01}$.  Now $(\tilde f_{02}\cup \tilde f_{21})(W^u_{012})$ is a curve in
$B_{0}\cup B_1$ of degree 2, and since $W^u_{0120}\cap B_{0,r}$ is 
connected, it
follows that $W^u_{0121}\cap B_{1,r}$ consists of two curves of degree 1.  By
Proposition 3.1, then, it follows that the complex variety 
$W^u_{0121}$ consists of
two irreducible components.  Now we map $W^u_{0121}\cap B_{1,r}$ under $\tilde
f_{12}$, which has degre 2.  By Proposition 3.4, $W^u_{01212}=\tilde
f_{12}W^u_{0121}$ lies inside $W^u_{012}$.  By $(**)$, $W^u_{01212}$ intersects
$W^s_{02}$.  Note that the arrangement of $W^u_{01212}$ corresponds 
to one of the
possibilities in Figure 7.  Another possiblity is given in the right 
hand of Figure 11.  This picture is mapped forward under $\tilde f_{20}$, to show one
possibility for $W^u_{012120}$ inside $B_{0,r}$.

\bigskip
\epsfxsize4.5in
\centerline{\epsfbox{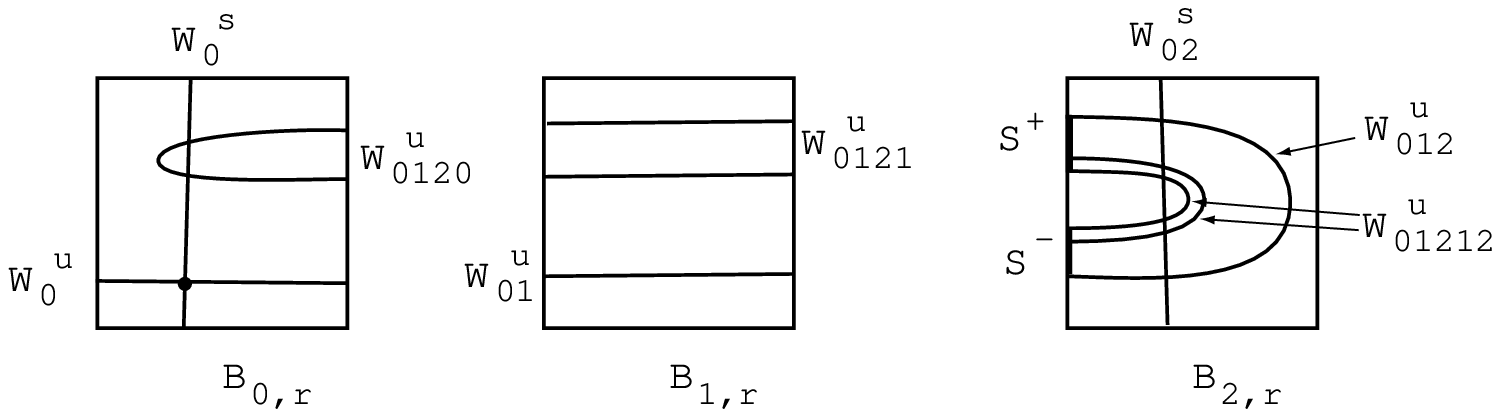}}
\centerline{Figure 10: Moving $W^u_0$ forward along the sequence 
01212 (case $b>0$)}

\bigskip
\epsfxsize3.2in
\centerline{\epsfbox{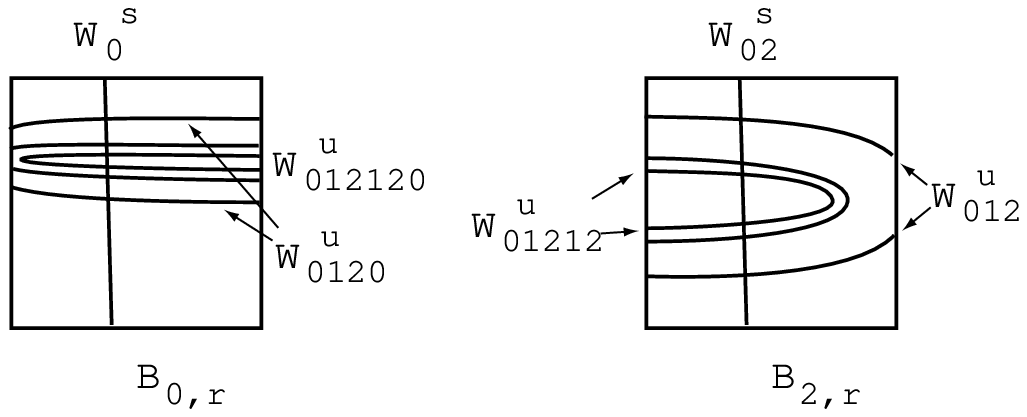}}
\centerline{Figure 11: Alternative to Figure 10}

Figure 12 deals with the orientation-reversing case and shows various unstable
pieces $W^u_I$ starting with $W^u_1$ through $p_1$ and moving forward along the
sequences $I=1200$, 1201, and 12012.  The construction of this picture 
was explained in
large part in the proof of Proposition 3.5, so we do not repeat it here.

\bigskip
\epsfxsize4.5in
\centerline{\epsfbox{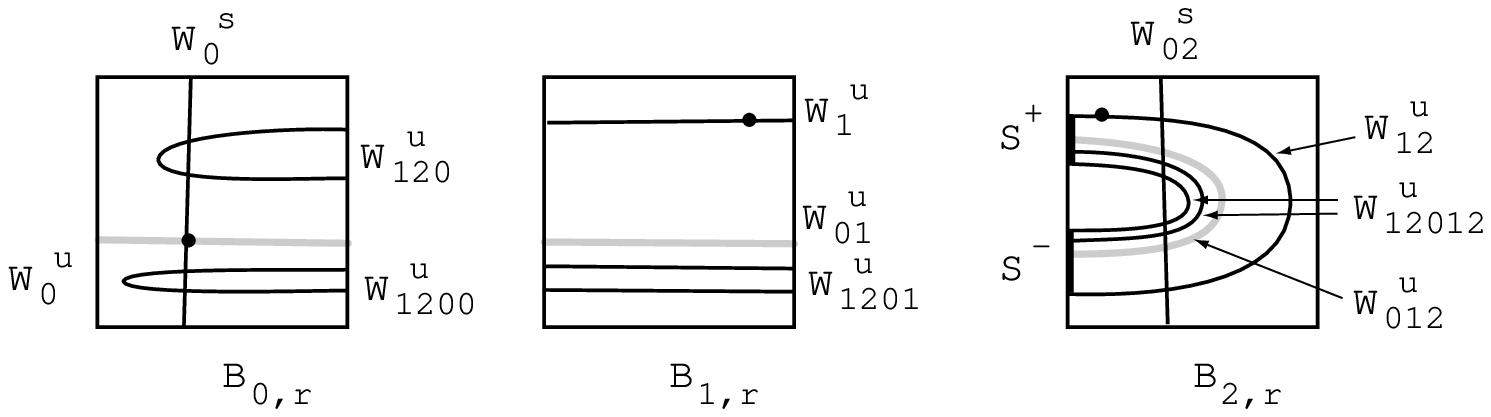}}
\centerline{Figure 12: Moving $W^u_1$ forward along the sequences 
1200, 1201, and
12012 (case
$b<0$)}\medskip
When $(**)$ holds, we use Figures 10 and 12 to define $S^\pm$ as the closed
subintervals of the left hand component of $\partial_vB_{2,r}$ which meet each
component of $\bar W^u_{012}\cup \bar W^u_{01212}$ if $b>0$ (resp.\ 
each component
of $\bar W^u_{12}\cup\bar W^u_{12012}$ if $b<0$).
\proclaim Proposition 3.6.  Suppose that $b>0$ and that 
(\dag) and $(**)$
hold. Let $I$ be an admissible sequence starting with 0 and ending 
with $k$, and
let $\Gamma$ be a connected component of $W^u_I$.  Then we have the following:
\vskip0pt If $k=0$, then: $\Gamma$ is disjoint from the component of
$B_{0,r}-W^u_0$ lying below $W^u_0$.  If $\delta(\Gamma)\ne1$,
then $\delta(\Gamma)=2$, and
$\Gamma$ intersects $W^s_0\cap B_{0,r}$, and $\bar\Gamma$ intersects 
the right hand
component of
$\partial_vB_{0,r}$ in two points.
\vskip0pt If $k=1$, then:  $\delta(\Gamma)=1$, and $\Gamma$ is 
disjoint from the
topmost and bottommost components of $B_{1,r}-(W^u_{01}\cup W^u_{0121})$.
\vskip0pt If $k=2$, then:  $\Gamma$ is disjoint from the innermost 
and outermost
components of $B_{2,r}-(W^u_{012}\cup W^u_{01212})$.  If
$\delta(\Gamma)\ne1$, then $\delta(\Gamma)=2$, and $\bar\Gamma$ intersects both
$S^+$ and
$S^-$.

\give Proof.  The proof proceeds by induction on the length of the 
sequence $I$.
First, the case $I=0$ is clear.  Now we suppose that the Proposition holds for
$I=I'i$.  We will show that if $(i,j)\in\cG$, then the Proposition holds for
$I=I'ij$ by considering five cases.

Case $ij=00$.  Since $f_{00}$ has sign $(+,\cdot)$, $f_{00}$ maps the 
component of
$B_{0,r}-W^u_0$ above $W^u_0$ to itself.  So $\tilde f_{00}\Gamma$ is 
disjoint from
the component of $B_{0,r}-W^u_0$ below $W^u_0$.  Now suppose 
$\delta(\Gamma)=2$.
$f_{00}$ maps $W^s_0$ into itself, and the sign of $f_{00}$ is $(\cdot,+)$, so
$\tilde f_{00}\Gamma$ intersects $W^s_0\cap B_{0,r}$, and 
$\bar\Gamma$ intersects
the right hand component of $\partial_vB_{0,r}$ in two points.

Case $ij=01$.  By Proposition 3.2, $f_{01}(B_{0,r})$ lies below 
$W^u_{0121}$.  On
the other hand $\Gamma$ is above $W^u_0$ and $f_{01}$ has sign $(+,\cdot)$, so
$\tilde f_{01}\Gamma$ is above $W^u_{01}$ in $B_{1,r}$.  It remains 
to show that
$\tilde f_{01}\Gamma$ consists of two components of degree 1.  For this, we may
assume that $\delta(\Gamma)=2$, and $\Gamma\cap W^s_0\cap B_{0,r}\ne\emptyset$.
Consider how $\gamma'=(\tilde f_{00}\cup \tilde 
f_{01})(B_{0,r}\cap\Gamma)$ maps
across
$B_{0,r}\cup B_{1,r}$: the left hand side of $\gamma'$ intersects $W^s_0\cap
B_{0,r}$ and the right hand side goes across the right hand boundary 
of $\partial_v
B_{1,r}$.  Thus $\gamma'\cap B_{1,r}$ consists of two curves.  By 
Proposition 3.1,
$\tilde g_{01}\Gamma\cap B_{0,1}=\gamma'\cap B_1$ consists of two 
disks of degree
one.

Case $ij=12$.  This is a direct consequence of Proposition 3.4.

Cases $ij=21$ and $ij=20$.  Let $\Gamma$ be as in case $k=2$.  We may 
assume that
$\delta(\Gamma)-2$.  Since $f_{20}$ and $f_{21}$ have sign 
$(\cdot,-)$, it follows
that  $\gamma':=(\tilde f_{20}\cup \tilde f_{21})\Gamma\cap 
(B_{0,r}\cup B_{1,r})$
is a 2-fold curve opening to the right.  Since $\bar\Gamma$ 
intersects both $S^+$
and $S^-$, we have $\Gamma \cap W^s_{02}\cap B_{2,r}\ne\emptyset$, 
and it follows
that $\gamma'\cap W^s_0\cap B_{0,r}\ne\emptyset$.  By Proposition 3.2, $\gamma'$
lies above $W^u_0$.  This finishes the case $ij=20$.

For the case $ij=21$, we observe that $\gamma'$ is a degree two real 
curve which
crosses $W^s_0\cap B_{0,r}$ and opens to the right.  Since 
$\bar\Gamma$ intersects
both $S^+$ and $S^-$, and $f_{21}$ has sign $(\cdot,-)$, it follows that
$\bar\gamma'$ intersects the right hand boundary of $\partial_vB_{1,r}$ in two
points.  Thus $\gamma'\cap B_{1,r}$ consists of two real curves 
crossing $B_{1,r}$
horizontally.  It follows from Proposition 3.1 that $\tilde 
f_{21}\Gamma$ consists
of two components of degree one.
\qed

In the following, we let $B^+_{0,r}$ denote the right-hand component of
$B_{0,r}-W^s_0$.
\proclaim Proposition 3.7.  Suppose that $b<0$, and that 
(\dag) and $(**)$
hold. Let $I$ be an admissible sequence starting with 1 and ending 
with $k$, and
let $\Gamma$ be a connected component of $W^u_I$.  Then we have the following:
\vskip0pt If $k=0$: $\Gamma$ is disjoint from the topmost and bottommost
components of $B^+_{0,r}-(W^u_{120}\cup W^u_{1200})$.  If 
$\delta(\Gamma)\ne1$, then
$\delta(\Gamma)=2$, and $\Gamma$ intersects $W^s_0\cap B_{0,r}$, and 
$\bar\Gamma$
intersects the right hand component of
$\partial_vB_{0,r}$ in two points.
\vskip0pt If $k=1$:  $\delta(\Gamma)=1$, and $\Gamma$ is disjoint from the
topmost and bottommost components of $B_{1,r}-(W^u_{1}\cup W^u_{1201})$.
\vskip0pt If $k=2$:  $\Gamma$ is disjoint from the innermost and outermost
components of $B_{2,r}-(W^u_{12}\cup W^u_{12012})$.  If
$\delta(\gamma)\ne1$, then $\delta(\Gamma)=2$, and $\bar\Gamma$ intersects both
$S^+$ and
$S^-$.

\give Proof.
This proof is analogous to the proof of Proposition 3.6; we omit the details.
\qed

\proclaim Proposition 3.8.  Suppose that (\dag) and $(**)$
hold.  Let $I$ be an admissible sequence of the form $I=0i_1\cdots 
i_n2$ if $b>0$
or $I=1i_1\cdots i_n2$ if $b<0$.  Then for each component $\Gamma$ of $W^u_I$,
$\#(W^s_{02}\cap \Gamma\cap B_{0,r})=\delta(\Gamma)$.  In particular, if the
intersection in the definition of $(**)$, is not tangential, then there is no tangency
between
$W^s_0$ and $W^u_I$.

\give Proof.  This follows from the case $k=2$ in Propositions 3.6 
and 3.7.  The
only case to consider is $\delta(\Gamma)=2$. Now if $\Gamma$ is not one
of the curves $W^u_I$ in condition $(**)$, $\Gamma\cap B_{2,r}$ is 
trapped between an
inner and an outer curve.  Since its closure intersects both
$S^+$ and
$S^-$, it must cross $W^s_{02}$ at least twice.  These two 
intersections account
for the total intersection number, and so these intersections must be simple
(nontangential), and there can be no further intersections.
\qed

\proclaim Proposition 3.9.  Suppose that (\dag) and $(**)$
hold.  Let $I$ be an admissible sequence of the form $I=0i_1\cdots 
i_n0$ if $b>0$
or $I=1i_1\cdots i_n0$ if $b<0$.  Then for each component $\Gamma$ of $W^u_I$,
$\#(W^s_0\cap \Gamma\cap B_{0,r})=\delta(\Gamma)$.  In particular, if the
intersection in the definition of $(**)$, is not tangential, then there is no tangency
between
$W^s_0$ and $W^u_I$.

\give Proof.  This follows by applying the map $\tilde f_{20}$, which 
has degree
one, to the result of Proposition 3.8.\qed

This allows us to characterize the mappings of maximal entropy.

\proclaim Theorem 3.9.  Suppose that (\dag) holds.  If the real map
$f_{a,b}$ has entropy equal to $\log 2$, then $(**)$ holds.    Conversely, if
$S\subset\{(a,b)\in\R^2:b\ne0\}$ is a connected set such that $(**)$  holds for all
$(a,b)\in S$, and if $f_{a_0,b_0}$ has entropy $\log 2$ for some  $(a_0,b_0)\in S$,
then $f_{a,b}$ has entropy $\log2$ for all $(a,b)\in S$.

\give Proof.  The proof will be  based on the following criterion from [BLS]:
$f_{a,b}$ has (maximal) entropy $\log2$ if and only if for all saddle 
points $p$ and
$q$, all (complex) intersection points of $W^s(p)\cap W^u(q)$ belong to $\R^2$.

We
suppose first that the entropy of $f_{a,b}$ is $\log2$.  If $b>0$ we take
$p=q=p_0$.  By (2.7), $\delta(W^s_{02})=1$ and $\delta(W^u_{01212})=4$.
If $b<0$, we take $p=p_0$ and $q=p_1$.  Again by (2.7), we have
$\delta(W^u_{12012})=4$. By (2.6) we have
$\#(W^s_{02}\cap\Gamma)=\delta(W^s_{02})\delta(\Gamma)=4$ with
$\Gamma=W^u_{01212}$ if $b>0$ and $\Gamma=W^u_{12012}$ if $b<0$.  By 
the criterion
above, all (complex) intersections between
$W^s_{02}$ and $\Gamma$ must belong to $\R^2$, so it follows that $(**)$ holds.

Now let us suppose that $(**)$ holds for all $(a,b)\in S$.  Consider the subset
$S_0$ of points $(a,b)\in S$ such that the entropy of $f_{a,b}$ is equal to
$\log 2$.  Since
$(a,b)\mapsto {\rm entropy}(f_{a,b})$ is continuous, it follows $S_0$ 
is a closed
subset of $S$.  Since
$S$ is connected, it suffices to show that $S_0$ is an open subset of $S$.
Let us fix a point $(a_0,b_0)\in S_0$.  By Proposition 2.5
there is an open set $U_0\subset\C^2$ such that $\bar B_{0,r}\cup\bar
B_{1,r}\cup\bar B_{2,r}\subset U_0$, and
$$U_0\cap f_{a,b}(B_0\cup B_1\cup B_2)\subset f_{a,b}(B_0\cup B_2)$$
$$U_0\cap  f_{a,b}(B_0\cup B_1\cup B_2)\subset f_{a,b}(B_1)$$
holds for $(a,b)=(a_0,b_0)$.  Thus it holds for $(a,b)$ in a small 
neighborhood of
$(a_0,b_0)$.   Thus we also have that 
$K_{a_0,b_0}\subset\R^2$ since
$f_{a_0,b_0}$ has maximal entropy.  By Proposition 2.1, then,
$K_{a_0,b_0}\subset B_{0,r}\cup B_{1,r}\cup B_{2,r}\subset U_0$.  Since
$(a,b)\mapsto K_{a,b}$ is upper semicontinuous, it follows that for $(a,b)$
sufficiently close to $(a_0,b_0)$ we have $K_{a,b}\subset U_0$, and thus
$f_{a,b}$ satisfies the hypotheses of Proposition 2.4.

Now we consider the case $b>0$; the argument for the case $b<0$ is 
similar and is
omitted.  Let
$q\in W^s(p_0)\cap W^u(p_0)$ be any point of intersection.  Replacing $q$ by
$f^{-m}q$ if necessary, we may assume that $q\in B_0$.  Let
$I$ denote the admissible sequence given by Proposition 2.4.  For $n$ 
sufficiently
large, we have
$f^nq\in W^s_{0}$, which is a neighborhood of $p_0$ inside $W^s(p_0)$.  Thus,
writing $I(n):=i_0i_1\cdots i_n$, we have
$f^nq\in W^u_{I(n)}$.  By Proposition 3.8, it
follows that
$W^s_{0}\cap W^u_{I(n)}\subset\R^2$.
Since $f^nq\in W^s_{0}\cap W^u_{I(n)}$, it follows that $q\in\R^2$.  Thus
$W^s(p_0)\cap W^u(p_0)\subset\R^2$, so that $f_{a,b}$ has entropy 
equal to $\log 2$.
\qed

\noindent{\bf Remark.}  There is an alternative approach to the ``Conversely'' part of this
Theorem.  Namely, we could use the arguments of this section to show that $W^s_{02}$ and
$W^u_I$ have certain trellis properties, and then we can apply the work of P. Collins [Co] to
conclude that the real map $f$ has entropy $\log 2$.

\section 4.  The Quadratic Horseshoe Locus

In this section we analyze the real, maximal entropy bifurcations in a neighborhood of
(2,0). 

\proclaim Lemma 4.1.  Suppose that $(a,b)\in\cA$, and suppose that
$\Delta/3\le\delta\le e^2-4-2\Delta$.  If we define $B_0'$ and $B_2'$ by
$$B_0':=\{|x+2|<\delta,|y|<e\}, \ \ B_2':=\{|x-2|<\delta,|y|<e\}$$ then $f$ induces crossed
mappings from $B_0'$ to itself and from 
$B_0'$ to $B_2'$. In particular, the sets $W^s_0$ and $W^s_{02}$ (as in (2.10)) are given by
$$W^s_0=\bigcap_{n\ge0}f^{-n}B_0', {\rm\ and\ }  W^s_{02}=\bigcap_{n\ge1}B_2'\cap
f^{-n}B_0'.$$

\give Proof.  Let us fix $\delta$ such that $\Delta/3\le \delta\le  e^2-4-2\Delta$ and set
$B_0':=\{|x+2|<\delta,|y|<e\}$.  By the upper bound on 
$\delta$, we have
$\{|4-x^2|<\Delta+\delta\}\subset\{|x|<e\}$.  We compute
$$\eqalign{f^{-1}B_0'\cap\{|y|<e\} &= \{|2+\pi_vf(x,y)|<\delta, |\pi_hf(x,y)|<e,|y|<e\}\cr
&\subset \{|2+a-x^2-by|<\delta,|x|<e,|y|<e\}\cr
&\subset\{|4-x^2|<|a-2|+|by|+\delta,|x|<e\}\subset\{|4-x^2|<\Delta+\delta,|x|<e\}
\cr &\subset\{|2-x|<\sqrt{4+\Delta+\delta}-2\}\cup
\{|2+x|<\sqrt{4+\Delta+\delta}-2\}\cr
&\subset\{|2-x|<{\Delta+\delta\over4}\}\cup\{|2+x|<{\Delta+\delta\over4}\}.\cr}$$ In the next
to last line we have removed the condition $|x|<e$ by the  upper bound condition on $\delta$. 
The last line uses the concavity of of the square root. By the lower bound on $\delta$, we
have $(\Delta+\delta)/4<\delta$,  so it follows that $f^{-1}\bar B_0'\cap\partial_v
B'_0=\emptyset$.

Next we consider a point $(x',y')\in f^{-1}(\partial_h B'_0)$.  By (?.5),
$$\eqalign{|y'|=|{1\over b}(a-y^2-x)| &> {e\over 
\Delta}(|y|^2-4-|a-2|-|x+2|)\cr &> {e\over\Delta}(e^2-4-\Delta-\delta).\cr}$$ This last
quantity is greater than $e$ by the upper bound on $\delta$, so
$(x',y')\notin\bar B_0'$.  Thus $f$ induces a crossed mapping from
$B_0'$ to itself.  The proof that $f$ induces a crossed mapping from $B_0'$ to
$B_2'$ is the same
\qed

\proclaim Corollary 4.2.  If $(a,b)\in\cA$, then $(*)$ holds.

\proclaim Proposition 4.3.  If $(a,b)\in\cA$, then the horizontal
disk $W^u_1$ has degree one.

\give Proof.  Let $\Gamma\in\cD^1_{h,r}(B_1)$ be any real, horizontal  disk.  Then by
Proposition 3.4, $\tilde f_{12}\Gamma\in\cD^2_{h,r}(B_2)$ is a real disk of degree two which
opens to the left. Applying $(\tilde  f_{20}\cup\tilde f_{21})$ to
$\tilde f_{12}\Gamma$, we obtain a disk $\Gamma'$ of degree two, which is horizontal in
$B_0\cup B_1$.  There can be at most one critical point for the projection $\pi_v:\Gamma'\to
B_0\cup B_1$, and if there is a critical point, it must be real, since its conjugate is also a
critical point.

Since the sign of $f_{20}\cup f_{21}$ is $(\cdot,-)$, $\Gamma'$ opens to the right.  By
Proposition 3.1, $\tilde f_{20}\tilde f_{12}\Gamma=\Gamma'\cap B_0$ defines a nonempty real
curve in $B_{0,r}$.  Thus, if there is a  critical point, then vertical projection
$\pi_v:\Gamma\cap(B_{0,r}\cup B_{1,r})\to(-d,c)$ has a critical  point.  Since $(*)$
holds, this critical point must belong to $B_{0,r}$, and by (\ddag), this point cannot
belong to $B_{1,r}$.  In particular, it follows that
$\pi_v$ has no critical point  in $\tilde f_{12}\tilde f_{12}\Gamma=\Gamma'\cap B_1$.  Thus
$\Gamma'\cap B_1$ consists of two components.  Since 
$p_1\in B_1\cap B_2$ is a fixed point, one of these components contains $p_1$, and  we denote
this component by $(\tilde f_{21}\tilde f_{12})^\#\Gamma$, which is a disk  of degree one.

Now if we choose $\Gamma$ to pass through $p_1$ such that its tangent  at $p_1$ is transverse
to $W^s(p_1)$, then .  It follows that
$(\tilde f_{12}\tilde f_{21})^{\#n}\Gamma$ is a sequence of horizontal disks of degree one,
passing through
$p_1$, which converge to $W^u_1$ as $n\to\infty$.  \qed

Now let us examine the case $b=0$.  The image of $f_{a,0}$ is the parabola
$$\Gamma:=f_{a,0}(\C^2)=\{x=a-y^2\}=\{(p(t),t):t\in\C\},$$
where $p(z)=a-z^2$.  Throughout our discussion, we assume that 
$|a-2|<\Delta$.  Thus
$a\notin D_0\cup D_1$, and so there are two holomorphic branches of
$p^{-1}(z)=\pm\sqrt{a-z}$ over
$D_0\cup D_1$.  For $j=0,1$, $\Gamma\cap B_j$ consists of two components
$\Gamma'_j$ and
$\Gamma''_j$, each of which is a smooth graph of a branch of $p^{-1}$.  We note
that $f_{a,0}$ is injective on each component $\Gamma'_j$ and $\Gamma''_j$,
$j=0,1$.  On the other hand, $\Gamma\cap B_2$ is connected, and $f_{a,0}$ is
two-to-one on $\Gamma\cap B_2$.

Let $p_0=(t_0,t_0)$ denote the fixed point which belongs to $B_0$.  (The following discussion
can be adapted to work with the other fixed point $p_1\in  B_1\cap B_2$, as well.)  Let
$\varphi_a:\C\to\C$ denote the linearizing coordinate such that 
$\varphi_a(0)=t_0$,
$\varphi_a'(0)=1$, and $p(\varphi_a(\zeta))=\varphi_a(\lambda\zeta)$, where
$\lambda:=p'(t_0)$.  If we write $\varphi=\varphi_a$, it follows that
$$\psi_{a,0}(\zeta):=(\varphi(\zeta),\varphi(\lambda^{-1}\zeta))$$ defines a mapping
$\psi_{a,0}:\C\to\Gamma$ which satisfies
$f_{a,0}\circ\psi(\zeta)=\psi(\lambda\zeta)$.

We wish to define the sets $W^u_I$ in the case $b=0$.  We let $W^u_0$ be the connected
component of
$\Gamma\cap B_0$ containing $p_0$;  $W^u_{01}$ is the connected component of
$\Gamma\cap B_1$ which intersects $W^u_0$; and $W^u_{012}=\Gamma\cap B_2$. As we  try to
consider longer $I$, we run into the difficulty that the mappings
$\tilde f_{ij}$ are not invertible.  To deal with this, we identify $W^u_I$ in terms of the
parametrization $\psi_{a,0}$ of $\Gamma$.  To do this, let
$\Omega_0\subset\C$ be the connected component of
$\psi^{-1}_{a,0}(W^u_0)=\varphi^{-1}(D_0)$ which contains the origin.  In general, we set
$$\Omega_I:=
\lambda^n\Omega_0\cap\varphi^{-n}D_{i_1}\cap\cdots\cap\varphi^{-n}(D_{i_n})
=\lambda^n\Omega_0\cap\psi_{a,0}^{-n}B_{i_1}\cap\cdots\cap\psi_{a,0}^{-n}(B_{i_n}),
\eqno(4.1)$$ where $I=0i_1\cdots i_n$ is an admissible sequence.  We then identify 
$W^u_I$ in terms of the map $\psi_{a,0}:\Omega_I\to W^u_I$.

The usefulness of the case $b=0$ is that it is the limit of the case $b\ne0$. When $b\ne0$, we
let
$\psi_{a,b}:\C\to W^u_{p_0}$ be the uniformization of $W^u(p_0)$,  normalized by the condition
$(\pi_v\circ\psi_{a,b})'(0)=1$.  In this case,
$(a,b,\zeta)\mapsto\psi_{a,b}(\zeta)$ is holomorphic, and we have
$$\lim_{b\to0}\psi_{a,b}=\psi_{a,0},\eqno(4.2)$$ with uniform convergence on compact subsets.
Restricting this to the image of $\Omega_I$, we have:
$$\lim_{b\to0}W^u_I(f_{a,b})=W^u_I(f_{a,0}),\eqno(4.3)$$ where the convergence is in the sense
of the Hausdorff topology.  Taking multiplicities of $W^{s/u}_I(f_{a,0})$ into account, the
convergence  also holds in the sense of currents.
\proclaim Lemma 4.4.  If $|a-2|<\Delta$, then for $I=01212$ and 
$I=12012$,
$\Omega_I$ consists of two connected components with disjoint closures.  If
$b\ne0$ is sufficiently small, then $W^u_I$ consists of two components.

\give Proof.  Since $p:D_0\to p(D_0)$ is a conformal equivalence, and
$p^{-1}D_0\subset D_0$, we may define a holomorphic map
$\lim_{n\to\infty}\lambda^np^{-n}:p(D_0)\to\C$  This is the inverse  of $\varphi$, and so
$\varphi:\lambda\Omega_0\to p(D_0)$ is univalent.  Thus 
$\Omega_{01}=\Omega_0\cap
\varphi^{-1}(D_1)$ is connected and relatively compact in $\Omega_0$.  Let $c_{01}$ be the
unique point of $\lambda\Omega_0$ such that 
$\varphi(c_{01})=0$.  It follows that $\varphi'(\lambda c_{01})=(p\circ\varphi(c_{01}))' = 
p'(0)\varphi'(c_0)=0$. Conversely, if $\zeta\in\lambda^2\Omega_0$, and if
$0=\varphi'(\zeta)= p'(\varphi(\lambda^{-1}\zeta)\varphi'(\lambda^{-1}\zeta)\lambda^{-1}$, 
then we must have $p'(\varphi(\lambda^{-1}\zeta))=0$ since $\varphi'\ne0$ on 
$\lambda\Omega_0$. It follows that $\zeta=\lambda c_{01}$, so $\lambda c_{01}$ is the  unique
critical point in $\lambda^2\Omega_0$.  It follows that
$\psi_{a,0}(\zeta)=(\varphi(\zeta),\varphi(\lambda^{-1}\zeta))$ has no critical point on
$\lambda^2\Omega_0$.  Since 
$\psi_{a,0}(\Omega_{012})=\Gamma\cap B_2$ is simply connected, it follows that
$\psi_{a,0}:\Omega_{012}\to\Gamma\cap B_2$ is univalent. By the argument above,
$\lambda^2c_{01}$ is the unique critical point for
$\psi_{a,0}$ in $\lambda^3\Omega_0$.  Now
$\psi_{a,0}(\lambda^2c_{01})=f_{a,0}(a,0)=(a-a^2,a)$, which does not belong to
$B_1$, since $\Re(a-a^2)<-c$.  It follows that $\psi_{a,0}$ is  unbranched on the closure of
$\lambda^3\Omega_0\cap\psi^{-1}_{a,0}(B_1)$.  Recall that
$f_{a,0}:W^u_{012}=\Gamma\cap B_2\to f_{a,0}(W^u_{012})$ is a mapping of degree two.  Thus
$W^u_{0121}$ is the component of $\Gamma\cap B_1$ which is  disjoint from
$W^u_{01}$, and $W^u_{0121}$ has multiplicity two.  It follows that
$\psi_{a,0}:\Omega_{0121}\to W6u_{0121}$ is a covering of degree two.  Since
$\psi_{a,0}$ is unbranched on the closure of
$\Omega_{0121}\subset\lambda^3\cap\psi^{-1}_{a,0}(B_1)$, it follows that
$\Omega_{0121}$ consists of two components with disjoint closures.

Let us move forward one more step: since
$f_{a,0}$ is injective on
$W^u_{0121}$, it follows that $\psi_{a,0}$ gives a conformal  equivalence between each
component of $\lambda\Omega_{0121}$ and $f_{a,0}W^u_{0121}$. Intersecting
$\lambda\Omega_{0121}$ with $\psi_{a,0}^{-1}(B_2\cap
f_{a,0}(W^u_{0121}))=\psi_{a,0}^{-1}(B_2\cap\Gamma)=\varphi^{-1}(D_2)$,  then
$\Omega_I$ consists of two components $\Omega_I'$ and
$\Omega_I''$ which have disjoint closures.  If $b\ne0$ is  sufficiently small, then
$\psi_{a,b}^{-1}(W^u_I)$ will be close to $\Omega_I$.  Thus it (as  well as $W^u_I$) has two
components.
\qed

Now we pass from unstable manifolds to stable manifolds.  The  vertical complex line through
the fixed point $p_0$ is mapped to $p_0$ under $f_{a,0}$.  If we write
$p_0=(t_0,t_0)$, then
$$W^s_0=\{(x,y):x=t_0,|y|<e\},{\rm\ and\ } W^s_{02}=\{(x,y):x=t_0',|y|<e\},$$ where
$t_0'\in\C$ is the solution to $p(t_0')=t_0$ such that $t_0'\ne t_0$. If $(\ddag)$ holds, then
$$W^s_{02}\cap\Gamma=\{(\zeta,\pm\sqrt{a-\zeta}):\zeta=t'_0\}.\eqno(4.4)$$ This intersection
consists of two distinct points unless $t'_0=a$,  which happens exactly when $a=2$.  We can
work our way backwards, taking successive  preimages, to define $W^s_J(f_{a,0})$ for an
admissible sequence $J$.  As in the  case of unstable manifolds, we have
$$\lim_{b\to0}W^s_J(f_{a,b})=W^s_J(f_{a,0}).\eqno(4.5)$$
\proclaim Proposition 4.5.  Suppose that $(a,b)\in\cA$ and
$|a-2|\ge(e+\Delta)|b|$.  Then for $I=01212$ and $I=12012$,
$W^s_{02}$ intersects $W^u_I$ in four distinct points, and thus the  intersection is not
tangential.

\give Proof.   We begin by noting
$$W^u_{I}\subset B_2\cap fB_0\subset
\{|b^{-1}(a-x-y^2)|<\delta,|y|<e\}.$$  If we set $\delta=\Delta/3$,  then by Lemma 3.5 we have

$$W^s_{02}\subset B_2'\subset\{|x-2|<{\Delta+\delta\over4},|y|<e\}.$$ Thus
$$\eqalign{W^s_{02}\cap W^u_{I} &\subset\{|x-2|<{\Delta+\delta\over4},|a-x-y^2|<\delta|b|\}\cr
&\subset\{|x-2|<{\Delta+\delta\over4},|a-2-y^2|<|b|\delta+|x-2|\}\cr
&\subset\{|a-2-y^2|<|b|\delta + {\Delta+\delta\over4}\}\cr
&=:U_{a,b}.\cr}$$
The set $U_{a,b}$ is symmetric with respect to $y\mapsto-y$ and is seen to be
disconnected if (and only if) it does not contain $y=0$.  This occurs exactly
when
$|a-2|\ge|b|\delta +
{\Delta+\delta\over4}$.  Now we recall that $\delta=\Delta/3$ and 
substitute the
condition
$(\ddag)$, which gives
$|a-2|\ge |b|\Delta/3 +\Delta/3\ge |b|\Delta/3 +(|a-2|+e|b|)/3$, and this is
equivalent to $|a-2|\ge(\Delta+e)|b|$.

Now consider the case $b=0$.  By Lemma 4.4, $\Omega_I$
consists of components
$\Omega_I'$ and
$\Omega_I''$.  Since $a\ne2$, the intersection (4.4) contains two
points, which lie in different components of $U_{a,b}$.  Thus
$\psi_{a,0}(\Omega'_I)$ and $\psi_{a,0}(\Omega_I'')$ each intersect 
$W^s_{02}$ in
two points, which lie in different components of
$U_{a,0}$.  If
$b\ne0$,
$|a-2|\ge(e+\Delta)|b|$, then $W^u_{I}$ consists of two components
${(W^u_I)}'=\psi_{a,b}(\Omega'_I(a,b))$ and
${(W^u_I)}''=\psi_{a,b}(\Omega''_I(a,b))$.  Further,  the set
$U_{a,b}$ continues to be disconnected, and by (4.3) the each component of
$U_{a,b}$ will continue to contain a point of
$W^s_{02}\cap {W^u_I}'$.  Since $\delta(W^u_I)'=2$ and  $W^s_{02}\cap 
{(W^u_I)}'$
contains two distinct points, the intersection is not tangential.  A similar
argument for $(W^u_I)''\cap W^s_{02}$ shows that $W^s_{02}\cap W^u_I$ has no
tangency.
\qed
Let us define $$\cD:=\{(a,b)\in\C^2: |a-2| < .237186,  |b|<.08205\}$$
$$T_{I}:=\{(a,b)\in\cD:W^s_{02}{\rm\ intersects\ } W^u_{I}{\rm\ 
tangentially}\}.$$
In the definition of $T_I$, we interpret the case $b=0$ as follows.
By \S1, we know that $T_I\cap\{b\ne0\}$ is a complex subvariety of
$\cD-\{b\ne0\}$.  By (4.5) and (4.2), we have that 
$(T_I\cap\{b\ne0\})\cup(2,0)$ is
the closure of $T_I\cap\{b\ne0\}$ in $\cD$.  With this interpretation,  $T_I$ is a
complex subvariety of $\cD$.

\proclaim Proposition 4.6.   For  $I=01212$
and $12012$, $T_I$ is a complex subvariety of $\cD$ with the 
following properties:
\vskip0pt (i) The projection
$\pi_h:T_I\to\{|b|<.08205\}$ is a proper mapping of degree two.
\vskip0pt (ii) $T_I$ is locally reducible at $(2,0)$.
\vskip0pt (iii)  There are real analytic functions
$\kappa^\pm_I:[-.08205,.08205]\to
\R$ with $\kappa_I^-(t)<\kappa_I^+(t)$ for $t>0$ such that
$T_I\cap\R^2$ is the union of the graphs of $\kappa^+_I$ and $\kappa^-_I$.

\give Proof.  Note that with our values of $e$ and $\Delta$, (\ddag) holds for
$(a,b)\in\cD$ whenever $b\ne0$.  Further, the condition
$|a-2|\ge(e+\Delta)|b|$ holds for $(a,b)\in\partial_v\cD$.  By Proposition 4.5,
then, $\bar T_I\cap \partial_v\cD=\emptyset$.  Thus $\pi_h$ is a 
proper mapping.
To determine the multiplicity of $\pi_h$, it suffices to determine 
the multiplicity
at
$b=0$.  If $b=0$, then the only tangency occurs at $a=2$.  Now
$W^u_{I}=\Gamma\cap B_2$, with multiplicity two, so in case $a=0$,
$W^s_{02}$ makes a tangential intersection with each component of
$W^u_I$.  It follows that $T_{I}\cap\{b=0\}=\{(2,0)\}$, with multiplicity two.
Thus $\pi_h$ has multiplicity two.

For (ii), let $b=0$.  By Lemma 4.4, $\Omega_I$ consists of components 
$\Omega_I'$
and $\Omega_I''$ which have disjoint closures.  Thus, for $b\ne 0$ 
small, there are
domains $\Omega_I'(a,b)$ and $\Omega''_I(a,b)$ which are mapped under 
$\psi_{a,b}$
to the two components of $W^u_I$.  Thus for $|b|<r_0$ small, we may split $T_I$
into $T_I'=\{(a,b)\in\cD:|b|<r_0, W^s_{02}{\rm\ intersects\ }W^u_I(a,b)'{\rm\
tangentially}\}$, and a similar set $T_I''$ for $W^u_I(a,b)''$.

Now we consider the projection 
$\pi_h:T_I\cap\R^2\to(-.08205,.08205)$.  This is a
proper mapping of degree two.  Consider a point $(a,b)\in 
T_I\cap\R^2$ with $b<0$
and suppose that $I=12012$.  We may repeat the argument of Proposition 3.5
to conclude that $W^u_{12012}$ consists of two curves in $B_{2,r}$ 
which open to
the left.  By $\gamma'$ and $\gamma''$ we denote the components of 
$W^u_{12012}$
such that
$\gamma'\cap B_{2,r}$ forms the inner curve, and $\gamma''\cap 
B_{2,r}$ forms the
outer curve.

Let us note at the outset that $\delta(\gamma')=\delta(\gamma'')=2$, and so
$\#(W^s_{02}\cap\gamma')=\#(W^s_{02}\cap\gamma'')=2$.  If there is a tangency
between $\gamma'$ and
$W^s_{02}$, then the tangency must be real.  For otherwise, if there 
were a point
of tangency
$q\in B_2-B_{2,r}$, the complex conjugate $\bar q$ would also be a point of
tangency, so the total intersection of
$\gamma'$ and $W^s_{02}$ in $B_2$ would be at least four.

Now suppose that the outer curve $\gamma''$ is tangential to
$W^s_{02}$.  Then this point of tangency must have order two, and can 
be the only
intersection with
$W^s_{02}$ since the total intersection satisfies $\#(W^s_{02}\cap
W^u_{12012})=2$.  Since $\gamma''\cap B_{2,r}$ opens to the left, it 
follows that
$\gamma''\cap B_{2,r}$ must lie to the left of
$W^s_{02}$.  Thus $\gamma'$ cannot intersect
$W^s_{02}\cap B_{2,r}$.  Thus there can be no tangency between the 
complex disks
$W^s_{02}$ and
$\gamma'$.

Thus in the case $b\ne0$, with $a$ and $b$ both real, there cannot be
tangencies (necessarily real) between both components of $W^u_I$ and 
$W^s_{02}$.
In other words, if
$(a,b)\in T'_I\cap\R^2$,
$b\ne0$, then
$(a,b)\notin T''_I\cap\R^2$.  This gives a splitting of $T_I$ into two
components in a neighborhood of $\pi_h^{-1}(-.08205,.08205)$.  Since
$\pi_h$ has degree one on $T'_I\cap\R^2$ and $T''_I\cap\R^2$ these sets
are given as the graphs of real analytic functions.
\qed
Let us set
$$\kappa(t):=\max(\kappa^+_{01212}(t),\kappa^-_{12012}(t)).$$
\proclaim Corollary 4.7.  $\{(a,b)\in\cD\cap\R^2:b\ne0,(**){\rm\ holds}\}
=\{(a,b)\in\cD\cap\R^2: b\ne0,a\ge\kappa(b)\}$.

\give Proof.  We consider only the case $b>0$; the other case is similar.  For
$I=01212$, set
$T_I^\pm:=\{a=\kappa^\pm_{01212}(b)\}$.  Thus
$T_I\cap\R^2=T_I^+\cup T^-_I$.    As was noted in the proof of Proposition~7.7,
$T^-_I$ is the set of parameters for which one component of $W^u_I$ 
is tangent to
$W^s_{02}$, and the other component is disjoint from $W^s_{02}$. 
$T^+_I$ is the
set of parameters for which one component of $W^u_I$ is tangential to 
$W^s_{02}$,
and the other component intersects $W^s_{02}$ in two points. \qed

Let us write
$$\cE:=\{(a,b)\in\R^2: f_{a,b}{\rm\ has\ entropy }<\log2\}$$
$$\cH:=\{(a,b)\in\R^2: f_{a,b}{\rm\ is\ a\ real\ horseshoe}\}$$

\proclaim Theorem 4.8. 
$$\cH\cap\cD=\{(a,b)\in\cD\cap\R^2:a>\kappa(b), b\ne0\},$$
$$\cE\cap\cD=\{(a,b)\in\cD\cap\R^2:a<\kappa(b),b\ne0\}.$$

\give Proof.  By Theorem 3.9 and Corollary 4.7, the set of parameters
$(a,b)\in\cD\cap\R^2$ for which the entropy is $\log 2$ is exactly the set
$\{a\ge\kappa(b)\}$.  On the other hand, if $a>\kappa(b)$, then by 
Proposition 3.8
there is no tangency.  Since $f$ has maximal entropy, it follows from 
[BS1] that
$f$ is hyperbolic.  Now $\cD\cap\R^2\cap\{a>\kappa(b)\}$ is a connected set of
parameters for which $f_{a,b}$ is hyperbolic.  By Theorem 1.1, this 
set contains
parameters for which $f_{a,b}$ is a real horseshoe.  It follows, then, from the
structural stability of hyperbolic maps that all of these maps are 
horseshoes.\qed

\section 5. Generic Unfolding

In Theorem 5.2 we establish the ``generic unfolding'' statement in Theorem 2. 
Let us fix $I=01212$ or $I=12012$.  In \S4 we saw that for $(a,b)\in\cD$, $b\ne0$, the set
$W^u_I$ is disconnected and may be split into 
$$W^u_I(a,b)=W^u_I(a,b)'\cup W^u_I(a,b)''.\eqno(5.1)$$
Further we saw that if
$(a_0,b_0)\in\cD\cap\R^2\cap\partial\cH$, then one of these components, say $W^u_I(a_0,b_0)'$,
has a quadratic tangency with $W^s_{02}(a_0,b_0)$.  This splitting may be done for all
$(a,b)\in\cD\cap\R^2$ in such a way that we obtain a continuous family
$$\cD\cap \R^2\ni (a,b)\mapsto W^u_I(a,b)'.$$ 
The horizontal projection $\pi_h(x,y)=y$,
establishes a conformal equivalence
$$\pi_h:W^s_{02}(a,b)\to\{|y|<e\}.$$  For $(a,b)\in\cD$, $b\ne0$, we define the function
$$h(a,b)=\prod_{i\ne j}(\pi_h(p_i)-\pi_h(p_j))$$ 
where the $p_i$ and $p_j$ in the product range over  the four points of
intersection $W^s_{02}(a,b)\cap W^u_I(a,b)$.  Since
$\pi_h|_{W^s_{02}(a,b)}$ is invertible, we see that
$h(a,b)\ne0$ if and only if there are four distinct points of intersection.  Thus $h(a,b)\ne0$
means that the multiplicities of all four intersections are 1, and thus all four
intersections are transverse.  As in \S4 we may extend the definition of $h$ to the case
$b=0$, and we see that
$h$ is analytic in $\cD$.

\proclaim Theorem 5.1.  For $(a,b)\in\cD\cap\R^2\cap T_I$ with $b\ne0$, we have ${\partial
h\over\partial a}\ne0$.

\give Proof. If $b=0$, then by the discussion in \S4, we see that
$a\mapsto h(a,0)$ has a zero of order 2 at $a=2$, and $h(a,0)\ne0$ for
$\{0<|a-2|<.237186\}$.  

By Theorem 4.5, none of the tangencies $T_I$ occur on the vertical boundary of $\cD$.  Thus
$h\ne0$ there.  It follows that for each fixed value $|b_0|\le .08$, the function
$$\{|a-2|<.237186\}\ni
a\mapsto h(a,b_0)$$ is analytic and has exactly two zeros (counted with
multiplicity).  One zero corresponds to a point $(a',b_0)\in T'_I$
and one corresponds to $(a'',b_0)\in T_I''$.  We have seen that
$T'_I\cap\cD\cap\R^2\cap\{b\ne0\}$ is disjoint from $T''_I\cap\cD\cap\R^2\cap\{b\ne0\}$. 
Since the total multiplicity is 2, each of these zeros must be a simple zero.  In particular,
we conclude that
${\partial h\over \partial z}(a,b)\ne0$ for $(a,b)\in T_I\cap\cD\cap\R^2\cap\{b\ne0\}$. \qed

Let us discuss this situation further.  We will consider a sequence of holomorphic
coordinate changes
$(x',y')=(x'(x,y),y'(x,y))$ which in addition depend holomorphically on the parameter
$(a,b)$.  First, we may change coordinates so that
$W^s_{02}(a,b)=\{x=0\}$ since
$W^s_{02}(a,b)$ has degree one in
$B_2$.  Now let
us split
$W^u_I(a,b)$ as in (5.1).  We will show that we may introduce coordinates such that we have
$$W^s_{02}=\{x=0\}{\rm\ \ and\ \  }W^u_I(a,b)=\{x=c_0(a,b)+y^2\}.\eqno(5.2)$$  The {\it generic
unfolding} condition is that ${\partial c_0(a,b)/ \partial a}\ne0$
for $a=a_b$ (see [PT, page 35]).

Now let us fix $b_0\in(-.08,.08)$, $b_0\ne0$, and set $a_0=a_{b_0}$.  Thus we have
$$W^u_I(a_0,b_0)'=\{x=\sum_{j=2}^\infty c_j(y-y_0)^j\},$$
where $(0,y_0)$ is the point of tangential intersection, and so $c_0=c_1=0$.  The
coefficient $c_2$ is nonzero because the intersection is quadratic (see [BS1]).  Without
loss of generality we may assume that
$y_0=0$.  Now for
$(a,b)$ near
$(a_0,b_0)$, we have
$$W^u_I(a,b)'=\{x=c_0(a,b) + c_1(a,b)y+c_2(a,b)y^2+\dots\}.$$
Now since $c_2(a,b)\ne0$ and $c_0(a_0,b_0)=c_1(a_0,b_0)=0$, we may solve $\tilde y=\tilde
y(a,b)\sim -c_1/(2c_2)$ such that
$${\partial x\over\partial y} =c_1(a,b)+2c_2(a,b)\tilde y+\dots=0.$$
Replacing $y$ by $y-\tilde y$, we have
$$W^u_I(a,b)'=\{x=\tilde c_0(a,b)+\tilde c_2(a,b)y^2+\dots\}.$$
Finally, since $\tilde c_2\ne0$, we may change coordinates $y'=\sigma(a,b)y$ to obtain (5.2).

Now we consider the function $h(a,b)$ in the coordinates $(x,y)$.  We have $W^s_{02}(a,b)\cap
W^u_I(a,b)'=\{(0,\pm\sqrt{-\tilde c_0(a,b)}\}$.  Since $W^u_I(a,b)'\cap
W^u_I(a,b)''=\emptyset$, and $W^u_I(a,b)''$ has no tangency for
$(a,b)$ near $(a_0,b_0)$, we have
$$h(a,b)=-(\sqrt{-\tilde c_0(a,b)}+\sqrt{-\tilde c_0(a,b)})^2\alpha(a,b)=2\tilde
c_0(a,b)\alpha(a,b)$$
where $\alpha$ is a nonvanishing analytic function.  Since $\tilde c_0(a_0,b_0)=0$, we have
$${\partial h\over\partial a}(a,b_0)={\partial \tilde c_0\over\partial
a}(a,b_0)\cdot \alpha(a,b_0)$$ 
for $a=a_0$.  By Theorem 5.1, then, ${\partial\tilde c_0}(a_0,b_0)/\partial a\ne0$.  Thus we
have:

\proclaim Theorem 5.2.  $(a,b)\mapsto (W^s_{02}(a,b),W^u_I(a,b))$ is a generic unfolding of a
tangency at the parameter value 
$(a_0,b_0)$.

\bigskip

\centerline{\bf References}

%\item{[A]}  L. Ahlfors, {\sl Conformal Invariants}, McGraw-Hill, New York,
%1973.
%
\item{[BLS]}  E. Bedford, M. Lyubich, and J. Smillie, Polynomial
diffeomorphisms of
$\C^2$.  IV: The measure of maximal entropy and laminar
currents,  Invent.\ Math.\ 112, 77--125 (1993).

%\item{[BS1]}  E. Bedford and J. Smillie,  Polynomial diffeomorphisms of
%$\C^2$: Currents, equilibrium measure, and hyperbolicity, Invent. Math. 87,
%69--99 (1990).

%\item{[BS6]}  E. Bedford and J. Smillie,  Polynomial diffeomorphisms of
%$\C^2$.  VI: Connectivity of $J$,  Annals of Math.\ 148, 695--735 (1998).
%
%\item{[BS7]}  E. Bedford and J. Smillie,  Polynomial diffeomorphisms of
%$\C^2$.  VII: Hyperbolicity and external rays,  Ann.\ Sci.\ Ecole Norm.\ Sup.,
%32, 455--497 (1999).

\item{[BS8]}  E. Bedford and J. Smillie,  Polynomial diffeomorphisms of
$\C^2$.  VIII: Quasi-expansion.  American J. of Math., 124, 221--271, (2002).

\item{[BS1]} E. Bedford and J. Smillie, Real polynomial diffeomorphisms with
maximal entropy: Tangencies.  Annals of Math., to appear.

\item{[BS3]} E. Bedford and J. Smillie, A new approach to the analysis of analytic
diffeomorphisms, in preparation.

%\item{[BiMi]} E. Bierstone and P. Milman, Semianalytic and 
%subanalytic sets, IHES
%Publications Math\'ematiques, 67, 5--42, 1988.

%\item{[CJY]}  L. Carleson, P. Jones, and J.-C. Yoccoz, Julia and John,
%Bol.\ Soc.\ Bras.\ Mat., Vol.\ 25, 1--30 (1994).
%
%\item{[Ca]}  J. Carette, Liens entre la g\'eometrie et la dynamique des
%ensembles
%de Julia, Th\`ese, U. de Paris-Sud, Orsay, 1997.

%\item{[Ch]} E. Chirka, {\sl Analytic Sets},  Kluwer Academic Publishers, 1989.

\item{[Co]} P. Collins, Dynamics forced by surface trellises. {\sl Geometry and topology
in dynamics (Winston-Salem, NC, 1998/San Antonio, TX, 1999)}, 65--86, Contemp. Math.,
246, Amer. Math. Soc., Providence, RI, 1999. 

\item{[DN]} R. Devaney and Z. Nitecki, Shift automorphisms in the 
H\'enon family,
Comm.\ Math.\ Phys.\ 67, 137--48, 1979.

%\item{[F]}  W.H.J. Fuchs, {\sl Topics in the Theory of Functions of 
%One Complex
%Variable}, van Nostrand, 1967.

\item{[FG1]}  J-E Forn\ae ss and E. Gavosto,  Existence of generic homoclinic tangencies for
H\'enon mappings, J.\ Geom.\ Anal., 2:5 (1992), 429--444.

\item{[FG2]}  J-E Forn\ae ss and E. Gavosto,    Tangencies for real and complex H\'enon maps:
An analytic method, Experimental Math., Vol. 8 (1999), No. 3, 253--260.

%\item{[FS]} J-E Forn\ae ss and N. Sibony,  Complex H\'enon mappings in $ C\sp 2$ and
%Fatou-Bieberbach domains. Duke Math. J. 65 (1992), no. 2, 345--380.

%\item{[FM]} S. Friedland and J. Milnor, Dynamical properties of plane polynomial
%automorphisms.  Ergodic Theory Dynam. Systems 9 (1989), no. 1, 67--99.

\item{[HO]} J.H.\ Hubbard and R.\ Oberste-Vorth, H\'enon mappings in the
complex domain II: Projective and inductive limits of polynomials, in: {\sl
Real and Complex Dynamical Systems}, B. Branner and P. Hjorth, eds. 89--132
(1995).

\item{[KKY]} ÊKan, Ittai; Kocak, Huseyin; Yorke, James A. Antimonotonicity: concurrent
creation and annihilation of periodic orbits. Ann. of Math. (2) 136 (1992), no. 2, 219--252. 

\item{[MNTU]} S. Morosawa, Y. Nishimura, M. Taniguchi, and T. Ueda, {\sl
Holomorphic Dynamics}, Cambridge U. Press, 2000.

\item{[O]} R. Oberste-Vorth, Complex horseshoes and the dynamics of 
mappings of two
complex variables, Ph.D dissertation, Cornell Univ., Ithaca, NY, 1987.

\item{[MT]} J. Milnor and W. Thurston, On iterated maps of the interval. {\sl Dynamical
systems (College Park, MD, 1986--87)}, 465--563, Lecture Notes in Math., 1342, Springer,
Berlin, 1988.

\item{[PT]}  J. Palis and F. Takens, {\sl Hyperbolicity and sensitive chaotic dynamics
at homoclinic bifurcations. Fractal dimensions and infinitely many attractors.}
Cambridge Studies in Advanced Mathematics, 35. Cambridge University Press, Cambridge,
1993.

\bigskip
\rightline{Eric Bedford}
\rightline{Indiana University}
\rightline{Bloomington, IN 47405}
\bigskip
\rightline{John Smillie}
\rightline{Cornell University}
\rightline{Ithaca, NY 14853}

\bye